\documentclass[final,a4paper,reqno]{elsarticle}
\usepackage{amssymb,amsmath,times,bbm,epsfig,float,fancyhdr}
\usepackage{amsmath}
\usepackage{wrapfig}
\usepackage{subfigure}
\usepackage{slashbox}

\graphicspath{{figures/}}

\usepackage{algorithm}
\usepackage{algorithmic}
\numberwithin{equation}{section}

\newcommand{\norm}[1]{\ensuremath{\left\|#1\right\|}}
\newcommand{\abs}[1]{\ensuremath{\left|#1\right|}}
\newcommand{\bu}{\ensuremath{\mathbf{u}}}

\newcommand{\mP}{\ensuremath{\mathcal{P}}}
\newcommand{\Grad}{\ensuremath{\mathrm{\nabla}}}

\newcommand{\Iona}{I_{1,\mathrm{ion}}}

\newcommand{\Ionb}{I_{2,\mathrm{ion}}}

\def\u {{\mathbf{u}}}

\newcommand{\T}{\mathcal{T}}

\newcommand{\br}{\boldsymbol{r}}

\newcommand{\R}{{\mathbb R}}

\def\grad{{\,\nabla }}

\def\Bleft{{\biggl[\hspace*{-3pt}\biggl[}}
\def\Bright{{\biggr]\hspace*{-3pt}\biggr]}}

\def\delt{{ {\scriptstyle \Delta}t }}

\def\D{{D}}

\newcommand{\Om}{\ensuremath{\Omega}}
\newcommand{\pt}{\ensuremath{\partial_t}}

\newcommand{\x}{\ensuremath{\mathbf{x}}}

\renewcommand{\u}{\mathbf{u}}

\newcommand{\Iap}{I_{\mathrm{app}}}
\newcommand{\Ion}{I_{\mathrm{ion}}}
\newcommand{\Iaph}{I_{\mathrm{app},h}}

\newcommand{\dx}{\ensuremath{\, dx}}

\newcommand{\dt}{\ensuremath{\, dt}}
\newcommand{\ds}{\ensuremath{\, ds}}

\newtheorem{thm}{Theorem}[section]
\newtheorem{rem}{Remark}[section]
\newtheorem{lem}[thm]{Lemma}
\newtheorem{defi}{Definition}[section]
\newtheorem{prop}{Proposition}[section]
\newenvironment{proof}[1][Proof]{\noindent\textit{#1.} }{\hfill \rule{0.5em}{0.5em}}
\allowdisplaybreaks

\begin{document}

\begin{frontmatter}
\title{A Virtual Element Method for a Nonlocal FitzHugh-Nagumo Model of Cardiac Electrophysiology}

\author[1]{Ver\'onica Anaya}
\ead{vanaya@ubiobio.cl}
\author[2]{Mostafa Bendahmane}
\ead{mostafa.bendahmane@u-bordeaux.fr}
\author[1,3]{David Mora}
\ead{dmora@ubiobio.cl}
\author[3,4]{Mauricio Sep\'ulveda}
\ead{mauricio@ing-mat.udec.cl}
\address[1]{GIMNAP, Departamento de Matem\'atica, Universidad del B\'io-B\'io,
Concepci\'on, Chile.}
\address[2]{Institut de Math\'ematiques de Bordeaux,
Universit\'e de Bordeaux, Talence, France.}
\address[3]{Centro de Investigaci\'on en Ingenier\'ia Matem\'atica
(CI$^2$MA), Universidad de Concepci\'on,
Concepci\'on, Chile.}
\address[4]{Departamento de Ingenier\'ia Matem\'atica,
Universidad de Concepci\'on, Concepci\'on, Chile.}

\begin{abstract}
We present a Virtual Element Method
(VEM) for a nonlocal reaction-diffusion system of the cardiac electric field.
To this system, we analyze an $H^1(\Om)$-conforming discretization
by means of VEM which can make use of general polygonal meshes.
Under standard assumptions on the computational
domain, we establish the convergence of the discrete
solution by considering a series of a priori estimates
and by using a general $L^p$ compactness criterion.
Moreover, we obtain optimal order space-time error estimates
in the $L^2$ norm.
Finally, we report some numerical tests supporting
the theoretical results.

\end{abstract}

\date{\today}
\end{frontmatter}

\section{Introduction}\label{sec:intro}
Reaction-diffusion systems appear in models of different areas such as medicine,
engineering, biology, physics, etc. The study of this kind of models has attracted too
much attention for many years, systems with different types of diffusion, for example: 
constant, nonlocal, cross. Mathematical models related with electrical activity 
in the heart (cardiac tissue) are becoming a powerful tools to study and understand
many types of heart disease, as for example irregular heart rhythm. 

The reaction-diffusion system of FitzHugh-Nagumo type \cite{FitzHugh61,nagumo}
is one of the most relevant and well-known generic model
in physiology which describes complex wave phenomena in
excitable or oscillatory media.
This model is a reaction-diffusion system which is a simplification of
the famous Hodgkin-Huxley model, which has been used to describe the propagation of
the electrical potencial in cardiac tissue \cite{Hastings75,Peskin:Book}.
The FitzHugh-Nagumo reaction-diffusion system consists
of one nonlinear parabolic partial differential equation (PDE)
which describes the dynamic of the membrane potential,
coupled with an ordinary differential equation which
models the ionic currents associated with the reaction term.
The main difficulties associated to solve this system are:
the coupling of the equations, through a nonlinear term and 
the regularity of the solution of the system is low.

In this paper, we analyze a {\it Virtual Element Method}
for a nonlinear parabolic problem arising in cardiac models (electrophysiology) with nonlocal diffusion
(see system \eqref{eq:main2} below). In our study, the self-diffusion coefficient is assumed depending on the total of electrical potential
in the heart.
The Virtual Element Method (VEM), recently introduced in
\cite{Beirao2,BBMR2014}, is a generalization of the
Finite Element Method which is characterized by the capability
of dealing with very general polygonal/polyhedral meshes. 
In recent years, the interest in numerical methods that can make
use of general polygonal/polyhedral meshes for the numerical
solution of partial differential equations has undergone a significant growth;
this because of the high flexibility that this kind of meshes allow
in the treatment of complex geometries. Among the large number of
papers on this subject, we cite as a minimal
sample~\cite{BLMbook2014,CGH14,DPECMAME2015,ST04,TPPM10}.

Although the VEM  is very recent, it has been applied to a large
number of  problems; for instance, VEM for Stokes, Brinkman,
Cahn-Hilliard, plates bending, advection-diffusion,
Helmholtz, parabolic, and hyperbolic
problems have been introduced in
\cite{ABMVsinum14,ABSV2016,BLV-M2AN,BMR2016,BM12,BBBPS2016,BGS17,CG16,CGS17,CMS2016,PPR15,vacca1,V-m3as18,vacca2},
VEM for spectral problems in \cite{BMRR,GVXX,MRR2015,MRV},
VEM for linear and non-linear
elasticity in \cite{ABLS2017a,BBM,BLM2015,Paulino-VEM,WRR2016}, 
whereas a posteriori error analysis have been developed
in \cite{BMm2as,BeBo2017,CGPS,MRR2}.

Over the past years, some papers related to numerical tools for
solving this model and its variations have appeared.
For example, in \cite{ChFPm2an13}
a continuous in space and discontinuous in time
Galerkin method of arbitrary order has been developed,
under minimal regularity assumptions, space-time error estimates
are established in the natural norms.
In \cite{Jackson92} some 
estimates in the $L^2$ norm for semi-discrete
Galerkin approximations for the FitzHugh-Nagumo model are derived.
A finite difference method has been presented in \cite{Barkley}, 
Chebyshev multidomain method has been presented in \cite{Olmos},
fully space-time adaptive multiresolution methods based on the
finite volume method and Barkley’s method for simulating the
complex dynamics of waves in excitable media in \cite{BRS}.
Finally, in \cite{Thomee:Book} has been presented
other methods related to the numerical analysis of general
semilinear parabolic PDE.

Numerical methods to solve these kind of models
have limitations in the range of applicable meshes.
In particular, finite element methods
rely on triangular (simplicial) or quadrilateral meshes.
Moreover, the classical finite volume method has
some restriction on the admisible meshes (for instance,
orthogonality constraints).
However, in complex simulations like fluid-structure
interaction, phase change, medical applications, and many others,
the geometrical complexity of the domain is a relevant issue
when partial differential equations have to be solved
on a good quality mesh; hence, it can
be convenient to use more general polygonal/polyhedral meshes.
Thus, in the present contribution, we are going to introduce and analyze
a VEM which has the advantage of using general polygonal meshes
to solve a nonlinear parabolic FitzHugh--Nagumo system,
where the diffusion coefficient depends on a nonlocal quantity.
The study of nonlocal diffusion problems has received considerable attention
in recent years since they appear in important
physical and biological applications \cite{ABLS15,ABS15,Chipot2000,ChipotLovat}. 
There are models of the FitzHugh-Nagumo type that also
take into account the nonlocal diffusion 
phenomena, for example in \cite{Liu2015} is considered
a diffusive nonlocal term as fractional diffusion,
in \cite{Oshita2003} is taken a nonlocal reactive term.

The aim of this paper is to introduce and analyze a conforming $H^1(\Om)$-VEM
which applies to general polygonal meshes, for the two-dimensional
nonlocal reaction-diffusion FitzHugh-Nagumo equations.
We propose a space discretization by means of VEM, which is
based on the discrete space introduced in \cite{AABMR13}
for the linear reaction--diffusion equation.
We construct a proper $L^2$-projection operator,
that is used to approximate the bilinear form that appears
for the time derivative discretization, which is obtained by
a classical backward Euler method. We also use that projection
to discretize the nonlocal term presented in the system.
We prove that the fully discrete scheme is well posed
and using standard space and time translates together with
a priori estimates for the discrete
solution, it is established convergence of the discrete
scheme to the weak solution of the model.
Under rather mild assumptions on the polygonal meshes,
we establish optimal order space-time error estimates
in the $L^2$ norm.

The structure of the paper is organized as follows:
in Section \ref{sec:DP}, we give some preliminaries
and assumptions on the data. Moreover, we introduce the
concept of weak solution. In Section \ref{section-vem},
we propose the semi-discrete and fully-discrete
virtual element method. In Section~\ref{sec-exist-est},
we prove the existence and convergence of the discrete solution.
In Section \ref{sec-error}, we give error
estimates, and finally, in Section \ref{sec:numerical-results},
some numerical results.

\section{Model problem and weak solution}\label{sec:DP}   

Fix a final time $T>0$ and a bounded domain $\Om \subset \mathbb{R}^2$ with polygonal boundary
$\Sigma$ and outer unit normal vector $\boldsymbol{n}$.
For all  $(x,t) \in \Omega_{T} :=\Omega \times (0,T)$, $v= v(x,t)$
and $w = w(x,t)$ stand for the transmembrane potential and the gating variable, respectively. 
The governing equations of the nonlocal reaction-diffusion FitzHugh-Nagumo system are:
\begin{eqnarray}
\label{eq:main2}
\begin{cases}
\displaystyle \partial_{t} v - \D \Biggl( \int_{\Om} v(x,t)\dx\Biggl)\Delta v+ 
\Ion(v,w) = \Iap(x,t) &\text{$(x,t) \in \Omega_{T}$,} \\
\partial_{t}w - H(v,w) = 0 &\text{$(x,t) \in \Omega_{T}$,}  \\
\\
\displaystyle \D \Biggl( \int_{\Om} v(x,t)\dx\Biggl)\Grad v\cdot \boldsymbol{n} = 0 &\text{$(x,t) \in \Sigma_{T} :=  \Sigma \times (0,T)$,} \\
\\
v(x,0) = v_{0}(x) &\text{$x \in \Omega$,} \\
w(x,0) = w_{0}(x) &\text{$x \in \Omega$.} 
\end{cases}
\end{eqnarray}
Herein, $\Iap$ is the stimulus. In this work, the diffusion rate $\D> 0$ is supposed to
depend on the whole of the transmembrane potential in the domain rather than on the local
diffusion, i.e. the diffusion of the transmembrane potential is guided by the global state of
the potential in the medium. We assume that
$\D: \R \to \R$ is a continuous function satisfying the following:
 there exist constants $d_1,\, d_2>0$ such that
\begin{equation}
\label{Cond_dcap3}
   d_1 \le \D \quad \text{and}\quad \abs{\D(I_1)-\D(I_2)} \le d_2\abs{I_1-I_2}
\quad \text{for all} \quad I_1, I_2 \in \R.
\end{equation}
Now, we make some assumptions on the data of the nonlocal FitzHugh-Nagumo model.
For the ionic current $\Ion(v,w)$, we assume that it can be decomposed into
$\Iona(v)$ and $\Ionb(w)$, where $\Ion(v,w) =\Iona(v) +\Ionb(w)$.
We assume that $\Iona,\Ionb:\R\to \R$ and  $H:\R^2\to \R$ are continuous functions, and that there
exist constants $\alpha_1, \alpha_2, \alpha_3, \alpha_4 >0$ such that
\begin{equation}\begin{split}
    \label{eq:ion1}
    &a)\quad\frac{1}{\alpha_1}\abs{v}^{4}\le \abs{\Iona(v)v}
    \le \alpha_2\left(\abs{v}^{4}+1\right),\\
    &b)\quad\abs{\Ionb(w)}\leq \alpha_3(\abs{w}+1),\\
    &c)\quad \forall \,z,s\in \R \quad  (\Iona(z)-\Iona(s))(z-s) \geq -C_h |z-s|^2,\\
    &d)\quad \abs{H(v,w)}\leq \alpha_4(\abs{v}+\abs{w}+1).
\end{split}\end{equation}

It is well known that the above assumptions
are fulfilled if the functions are specified as follows:
\begin{eqnarray}
\label{eq:H}
H(v,w)=av - bw,
\end{eqnarray}
and
\begin{eqnarray}
\label{eq:ionplus}
\Ion(v,w) = -\lambda(w-v(1-v)(v-\theta)),
\end{eqnarray} 
where $a,b,\lambda,\theta$ are given parameters.

Next, we will use the following spaces: by $H^m(\Om)$,
we denote the usual Sobolev space of order $m$. 
Given $T>0$ and $1\leq p \leq
\infty$, $L^p(0,T;\mathbb{R})$ denotes the space of $L^p$ integrable
functions from the interval $[0,T]$ into $\mathbb{R}$. 
The weak solution to the model (\ref{eq:main2}) is defined as follows. 
\begin{defi}[Weak solution]
A weak solution to the system (\ref{eq:main2}) is a double function $(v, w)$ such that $v \in L^{2}(0,T;H^{1}(\Omega)) \cap L^{4}(\Omega_{T}) $,
$\partial_{t}v\in L^{2}(0,T;(H^{1}(\Omega)')) + L^{\frac{4}{3}}(\Omega_{T})$, $w \in C([0,T];L^{2}(\Omega)) $, and satisfying the following weak formulation  
\begin{equation}\begin{split}\label{eq:weak-sol-1}
\iint_{\Omega_{T}}{\partial_{t} v\,\varphi} + \int_0^T \D \Biggl( \int_{\Om} v(x,t)\dx\Biggl)
\int_{\Omega} { \nabla v \cdot \nabla \varphi} + \iint_{\Omega_{T}}{\Ion(v,w) \varphi}
&= \iint_{\Omega_{T}} {\Iap(x,t)\varphi},  \\
\iint_{\Omega_{T}} {\partial _t w \, \phi -\iint_{\Omega_{T}}  H(v,w) \phi} &=0,
\end{split}\end{equation}
for all $\varphi \in L^{2}(0,T;H^{1}(\Omega))\cap L^{4}(\Omega_{T})$ and $\phi \in C([0,T];L^{2}(\Omega))$.
\label{defi-ws}
\end{defi}

\begin{rem}
Note that, in view of the conditions stated in \eqref{eq:ion1}, we can easily check
that Definition \ref{defi-ws} makes sense. Furthermore, observe that
Definition \ref{defi-ws} implies $v \in C([0, T];L^2(\Om))$ (see \cite{Schoenbek78}).
\end{rem}

\setcounter{equation}{0}
\section{Virtual element scheme and main result}\label{section-vem}
In this section, first we recall the mesh construction and the
assumptions consider to introduce the discrete virtual element space.
Then, we present the virtual element approximation of the FitzHugh-Nagumo model.
In the sequel, the existence and uniqueness is proved.

\subsection{The VEM semi-discrete problem}\label{sec:femdef}

Let $\left\{{\mathcal T}_h\right\}_h$ be a sequence of decompositions of $\Omega$
into polygons $K$. Let $h_K$ denote the diameter of the element $K$
and $h$ the maximum of the diameters of all the elements of the mesh,
i.e., $h:=\max_{K\in{\mathcal T}_h}h_K$. In what follows,
we denote by $N_K$ the number of vertices of $K$,
by $e$ a generic edge of $\mathcal{T}_h $
and for all $e\in \partial K$, we define a unit normal vector $\boldsymbol{n}_K^e$ 
that points outside of $K$.

For the analysis, we will make the following
assumptions as in \cite{Beirao2,BMRR}:
there exists a positive real number $C_{{\mathcal T}}$ such that,
for every $h$ and every $K\in {\mathcal T}_h$,
\begin{itemize}
\item[{\bf A1}:] the ratio between the shortest edge
and the diameter $h_K$ of $K$ is larger than $C_{{\mathcal T}}$;
\item[{\bf A2}:] $K\in{\mathcal T}_h$ is star-shaped with
respect to every point of a  ball
of radius $C_{{\mathcal T}}h_K$.
\end{itemize}

For any subset $S\subseteq\R^2$ and nonnegative
integer $k$, we indicate by $\mathbb{P}_{k}(S)$ the space of
polynomials of degree up to $k$ defined on $S$.

Now, we consider a simple polygon $K$
(meaning open simply connected set whose boundary is a non-intersecting line
made of a finite number of straight line segments), and we start by
introducing a preliminary virtual element space. For all $K \in {\mathcal T}_h$,
the local space $V_{k|K}$ is defined by
$$V_{k|K}:=\{\varphi \in H^1(K)\cap C^0(K) : \varphi_{|e} \in \mathbb{P}_k(e) \; \forall e \in \partial K,
 \; \Delta \varphi \in \mathbb{P}_k(K) \}.$$
Now, we introduce the following set of linear operators
from $V_{k|K}$ into $\R$. For all $\varphi\in V_{k|K}$:
\begin{itemize}
\item $D_1$: The values of $\varphi$ at the vertices of $K$;
\item $D_2$: Values of $\varphi$ at $k-1$ distinct points
in $e$, for all $e\in\partial K$;
\item $D_3$: All moments $\displaystyle \int_{K}\varphi p\,dx$, for all $p\in \mathbb{P}_{k-2}(K)$.
\end{itemize}

Now, we split the bilinear form $a(\cdot,\cdot):=(\Grad \cdot , \Grad \cdot )_{0,\Om}$,
\begin{equation*}
a(v,\varphi):= \sum_{K \in {{\mathcal T}_h}} a^K(v,\varphi),
\qquad \forall v,\varphi \in H^1(\Om),
\end{equation*}
where
\begin{equation*}
a^K(v,\varphi):= \int_K\nabla v\cdot\nabla\varphi,
\qquad \forall v,\varphi \in H^1(\Om).
\end{equation*}

For the analysis we will introduce the following
broken seminorm:
$$|\varphi|_{1,h}:=\left( \sum_{K \in {{\mathcal T}_h}} |\varphi|^2_{1,K} \right)^{1/2}.$$

Let $\Pi_{K,k} : V_{k|K} \to  \mathbb{P}_k(K)$ be the projection operator defined by

\begin{equation*}
\begin{cases}
& \displaystyle a^K(\Pi_{K,k}v,q)=a^K(v,q) \qquad \forall q \in \mathbb{P}_k(K),\\
& P_0(\Pi_{K,k}v)=P_0 v,
\end{cases}
\end{equation*}
where $P_0$ can be taken as
\begin{equation*}
\begin{cases}
&\displaystyle   P_0v:=\frac{1}{N_K}\sum_{n=1}^{N_K}v(V_i) \qquad k=1, \\
& P_0v:=\frac{1}{|K|} \displaystyle \int_K v \dx \qquad k>1,
\end{cases}
\end{equation*}
with $V_i$ the vertices of $K$,  $1 \leq i \leq N_K$
where $N_K$ is the number of vertices in $K$.

Using an integration by parts, it is easy to check that,
for any $\varphi\in V_{k|K}$, the values of
the linear operators $D_1,D_2$ and $D_3$ given before are sufficient
in order to compute $\Pi_{K,k}$. As a consequence, the
projection operator $\Pi_{K,k}$ depends only on the
values of the operators $D_1,D_2$ and $D_3$.

Now, we introduce our virtual local space
$$W_{k|K}:=\left\{\varphi \in V_{k|K} : \int_K (\Pi_{K,k} \varphi)q \dx = \int_K \varphi q \dx \quad \forall q \in 
\mathbb{P}_k/\mathbb{P}_{k-2}(K) \right\},$$
where the symbol $\mathbb{P}_k/\mathbb{P}_{k-2}(K)$ denotes the polynomials
of degree $k$ living on $K$ that are $L^2$-orthogonal to all polynomials
of degree $k-2$ on $K$. We observe that, since $W_{k|K}\subset V_{k|K}$,
the operator  $\Pi_{K,k}$ is well defined on $W_{k|K}$ and computable
only on the basis of the values of the operators $D_1,D_2$ and $D_3$.

The global discrete space will be

$$W_h:=\{\varphi \in H^1(\Omega) : \varphi|_K \in W_{k|K}, \quad \forall K \in {\mathcal T}_h \}.$$
In agreement with the local choice of the degrees of freedom, in
$W_h$ we choose the following degrees of freedom:
\begin{itemize}
\item $DG_1$: The values of $\varphi$ at the vertices of $\mathcal{T}_h$;
\item $DG_2$: Values of $\varphi$ at $k-1$ distinct points
in $e$, for all $e\in\mathcal{T}_h$;
\item $DG_3$: All moments $\displaystyle \int_{K}\varphi p\,dx$, for all $p\in \mathbb{P}_{k-2}(K)$
on each element $K\in\mathcal{T}_h$.
\end{itemize}

On the other hand, let $S^K(\cdot, \cdot)$ and $S_0^K(\cdot, \cdot)$
be any symmetric positive definite bilinear
forms to be chosen as to satisfy
\begin{align}
&c_0 a^K(\varphi_h,\varphi_h)\leq S^K(\varphi_h,\varphi_h) \leq c_1 a^K(\varphi_h,\varphi_h)
\quad \forall \varphi_h \in V_{k|K} \quad {\text with} \quad \Pi_{K,k} \varphi_h=0,\label{Sestab1}\\
&\tilde{c}_0 (\varphi_h,\varphi_h)_{0,K}\leq S_0^K(\varphi_h,\varphi_h)
\leq \tilde{c}_1 (\varphi_h,\varphi_h)_{0,K} \quad \forall \varphi_h \in V_{k|K},\label{Sestab2}
\end{align}
for some positive constants $c_0,c_1,\tilde{c}_0$ and $\tilde{c}_1$ independent of $K$. 

We define the local discrete bilinear and trilinear forms:

$$a^K_h(\cdot,\cdot):W_h \times W_h \to \R, \qquad  m^K_h(\cdot,\cdot):W_h \times W_h \to \R,$$
$$b^K_h(\cdot,\cdot, \cdot ):W_h \times W_h\times W_h \to \R,
\quad c^K_h(\cdot,\cdot,\cdot):W_h \times W_h\times W_h \to \R,$$
as follow, for all $v_h, w_h, \varphi_h \in W_{k|K}$:

$$a_h^K(v_h,\varphi_h) := a^K(\Pi_{K,k} v_h,\Pi_{K,k} \varphi_h)+
S^K(v_h - \Pi_{K,k}v_h, \varphi_h - \Pi_{K,k} \varphi_h),$$

$$m_h^K(v_h,\varphi_h) := (\Pi^0_{K,k} v_h,\Pi^0_{K,k} \varphi_h)+
S_0^K(v_h - \Pi^0_{K,k}v_h, \varphi_h - \Pi^0_{K,k}\varphi_h),$$

$$b_h^K(v_h,w_h,\varphi_h) := \int_K I_{ion}(\Pi^0_{K,k} v_h,\Pi^0_{K,k} w_h)\Pi^0_{K,k} \varphi_h,$$

$$c_h^K(v_h,w_h,\varphi_h) := \int_K H(\Pi^0_{K,k} v_h,\Pi^0_{K,k} w_h)\Pi^0_{K,k} \varphi_h,$$
where $\Pi^0_{K,k}:W_{k|K} \to \mathbb{P}_k(K)$ is the
standard $L^2$-projection operator which is
computable on the basis of the degrees of freedom (see \cite{AABMR13,vacca2}).

We observe that for all $K \in {\mathcal T}_h$ it holds:
\begin{itemize}
\item $k$-consistency: for all $p \in \mathbb{P}_k(K)$ and for all $\varphi_h \in W_{k|K}$
\begin{equation}\label{consis}
\begin{split}
a^K_h(p,\varphi_h)=a^K(p,\varphi_h),\\
m^K_h(p,\varphi_h)=(p,\varphi_h)_{0,K}.
\end{split}
\end{equation}
\item stability: there exist four positive constants $\alpha', \alpha'', \beta', \beta''$, independent
of $h$, such that for all $\varphi_h \in W_{k|K}$
\begin{equation}\begin{split}\label{stability}
&\alpha' \; a^K(\varphi_h,\varphi_h) \leq a^K_h(\varphi_h,\varphi_h) \leq  \alpha''  \; a^K(\varphi_h,\varphi_h),\\
&\beta'  \; (\varphi_h,\varphi_h)_{0,K} \leq m^K_h(\varphi_h,\varphi_h) \leq  \beta''  \; (\varphi_h,\varphi_h)_{0,K}.
\end{split}\end{equation}
\end{itemize}

Then, we set for all $v_h, w_h, \varphi_h \in W_h$,
$$a_h(v_h, \varphi_h):= \sum_{K \in {{\mathcal T}_h}} a^K_h(v_h,\varphi_h),\quad m_h(v_h, \varphi_h):= \sum_{K \in {{\mathcal T}_h}} m^K_h(v_h,\varphi_h),$$ 

$$b_h(v_h,w_h, \varphi_h):= \sum_{K \in {{\mathcal T}_h}} b^K_h(v_h,w_h,\varphi_h), \quad c_h(v_h,w_h, \varphi_h):= \sum_{K \in {{\mathcal T}_h}} c^K_h(v_h,w_h,\varphi_h).$$ 

We discretize the nonlocal diffusion term
using the $L^2$-projection
as follows:
\begin{equation}\label{defJ}
J(v_h):=\int_{\Om}v_h=\displaystyle \sum_{K \in {\mathcal T}_h}
\displaystyle \int_K \Pi^0_{K,k} v_h, \qquad v_h \in W_h.
\end{equation}

For the right-hand side, we assume $\Iap(x,t) \in L^{2}(\Om_{T})$ and we set
$$I_{app,h}(t) =  \Pi^0_{k} \Iap(\cdot,t) \quad {\mathrm for\,\,a.e.} \quad  t \in (0,T),$$
where we have introduced
$\Pi_k^0$ as the following operator which is defined in $L^2$ by 
\begin{equation}
\label{proyectorg}
(\Pi_k^0g)|_K:=\Pi_{K,k}^0g\quad\text{ for all }K\in \mathcal{T}_h
\end{equation}
with $\Pi_{K,k}^0$ the $L^2(K)$-projection.

Now, we note that from the symmetry of $a_h(\cdot,\cdot)$
and $m_h(\cdot,\cdot)$ and the stability conditions stated
before imply the continuity of $a_h$ and $m_h$. In fact,
for all $v_h,\varphi_h\in W_h$:
\begin{equation}\label{continuity}
\begin{split}
\vert a_h(v_h,\varphi_h)\vert\le C\Vert v_h\Vert_{H^1(\Om)}\Vert\varphi_h\Vert_{H^1(\Om)},\\
\vert m_h(v_h,\varphi_h)\vert\le C\Vert v_h\Vert_{L^2(\Om)}\Vert\varphi_h\Vert_{L^2(\Om)}.
\end{split}
\end{equation}

The semidiscrete VEM formulation reads as follows.
For all $t>0$, find $v_h,w_h \in L^2(0,T;W_h)$ with $\partial_t v_{h}, \partial_t w_{h} \in L^2(0,T;W_h)$ such that
\begin{equation}\begin{cases}\label{eq:FEM1}
  &\displaystyle   m_h(\partial_t v_{h}(t),\varphi_h) +\D\left( J(v_h(t))\right)
  a_h \Bigl( v_h(t),\varphi_h\Bigr)
  + b_h (v_h(t),w_h(t),\varphi_h)
  = \left( I_{app,h}(t),\varphi_h \right)_{0,\Om}\\
&\displaystyle  m_h(\partial_t w_{h}(t),\phi_h) - c_h (v_h(t),w_h(t),\phi_h)=0,
\end{cases}\end{equation}
for all $\varphi_h,\phi_h\in W_h$. Additionally, we set  $v_h(0)=v_{h}^0$
and $w_h(0)=w_{h}^0$.
A classical backward Euler integration method is employed
for the time discretization of \eqref{eq:FEM1} with time step $\Delta t=T/N$.
This results in the following fully discrete method:
find  $v_h^n,w_h^n\in W_h$ such that

\begin{equation}\begin{cases}\label{wf-finite}
&\displaystyle m_h \left(\frac{v_{h}^n\!-\!v_{h}^{n-1}}{\Delta t},\varphi_h\right)
+\D \left( J(v_{h}^n) \right)a_h\Bigl(v_{h}^n,\varphi_h\Bigr) + b_h (v_h^n,w_h^n,\varphi_h)
=\left( I_{app,h}^n,\varphi_h \right)_{0,\Om}\\
&\displaystyle m_h\left( \frac{w_{h}^n\!-\!w_{h}^{n-1}}{\Delta t},\phi_h\right) - c_h (v_h^n,w_h^n,\phi_h)=0,
\end{cases}\end{equation}
for all $\varphi_h,\phi_h\in W_h$, for all $n\in \{1,\ldots,N\}$;
the initial condition takes the form $v_{h}^0, w_{h}^0$
and $I_{app,h}^n:=I_{app,h}(t_n)$ with $t_n:=n\Delta t$, for $n=0,\ldots,N$.

Our main result is the following theorem:
\begin{thm} \label{theo1}
Assume that \eqref{Cond_dcap3} and \eqref{eq:ion1} hold.
If $v_{0}(x)\in L^2(\Om)$, $w_{0}(x)\in L^2(\Om)$,  and $\Iap(x,t) \in L^{2}(\Om_{T})$, 
then the virtual element solution $\bu_h^n=\Bigl(v_h^n,w_h^n\Bigr)$, generated 
by \eqref{wf-finite}, converges 
along a subsequence to $\bu=(v,w)$ as $h\to 0$, 
where $\bu$ is a weak solution of \eqref{eq:main2}. 
Moreover, the weak solution is unique.
\end{thm}

In the next section, we prove Theorem \ref{theo1} by establishing
the convergence of the virtual element solution
$\Bigl(v_h^n,w_h^n\Bigr)$, based on a priori estimates and the compactness
method. Moreover, we provide error estimates in Section \ref{sec-error}.

\setcounter{equation}{0}
\section{Existence of solution for the virtual element scheme}\label{sec-exist-est}

The existence result for the virtual element
scheme is given in the following proposition.

\begin{prop}\label{prop-fem}
Assume that \eqref{Cond_dcap3} and  \eqref{eq:ion1} hold.
Then, the problem \eqref{wf-finite} admits a discrete solution
$\bu_h^n=\Bigl(v_h^n,w_h^n\Bigr)$.
\end{prop}
\begin{proof}
The existence of $\bu_h^n$ is shown by induction on $n=0,\dots,N$. For
$n=0$, solution is given by $\bu_h^0=(v_h(0),w_h(0))=(v_{h}^0,w_{h}^0)$. Assume that
$\bu_h^{n-1}$ exists. Choose $\Bleft\cdot\,,\,\cdot\Bright$
as  the scalar product on $H^{1}(\Om)\times L^2(\Om)$. We are looking for a
solution $\bu_h^{n}$ to $\Bleft L(\bu_h^{n}),\Phi_h\Bright=0$, where the
operator $L:W_h \times W_h \to W_h \times W_h$ is given by
\begin{align*}
\Bleft L\bigl(\bu_h^{n}\bigr),\Phi_{h}\Bright=&m_h \left(\frac{v_{h}^n\!-\!v_{h}^{n-1}}{\Delta t},\varphi_h\right)
+\D \left( J(v_{h}^n) \right)a_h\Bigl(v_{h}^n,\varphi_h\Bigr)+b_h(v_h^n,w_h^n,\varphi_n)- \left ( I_{app,h}(t_n),\varphi_h \right)_{0,\Om}\\
&\displaystyle \quad +m_h\left( \frac{w_{h}^n\!-\!w_{h}^{n-1}}{\Delta t},\phi_h\right)- c_h(v_{h}^n,w_{h}^n,\phi_h),
\end{align*}
for all $\Phi_h:=(\varphi_h,\phi_h)\in W_h \times W_h$.
Note that the continuity of the operator $L$ is a consequence of the continuity of $m_h$, $a_h$ $b_h$ and $c_h$. Moreover, the following bound holds from the discrete H\"{o}lder
and Sobolev inequalities (recall that $H^1(\Om)\subset L^q(\Om)$ for all $1\leq q \leq 6$):
\begin{align*}
\Bleft L\bigl(\u^{n}_h\bigr),\Phi_{h}\Bright\le C (\norm{v^{n}_h}_{H^{1}(\Om)}+\norm{w^{n}_h}_{L^2(\Om)}+1)
(\norm{\varphi_h}_{H^{1}(\Om)}+\norm{\phi_h}_{L^2(\Om)}),
\end{align*}
for all $\bu_h$ and $\Phi_{h}$ in $W_h \times W_h$.
Moreover, from \eqref{eq:ion1} and Young inequality, we get
\begin{equation*}\label{eq:disc-exist:3}
\begin{split}
\Bleft L(\bu^{n}_h),\bu^{n}_h\Bright\geq C(\norm{v^{n}_h}_{H^{1}(\Om)}^2
+\norm{w^{n}_h}_{L^2(\Om)}^2)+C'
\end{split}
\end{equation*}
for some constant $C,C'>0$.
Finally, we conclude that
$\Bleft L(\bu^{n}_h) \,,\,\bu^{n}_h\Bright\geq 0$ for
$\|\bu^{n}_h\|^2:=\norm{v_h^n}_{H^{1}(\Om)}^2+\norm{w_h^n}_{L^2(\Om)}^2$
sufficiently large. The existence of
$\bu^{n}_h$ follows by the standard Brouwer fixed point
argument (see \cite[Lemma 4.3]{Lions:Book69_Fr}).
\end{proof}

\subsection{A priori estimates}\label{sec:basic-apriori}

In this section, we establish several a priori (discrete energy) estimates 
for the virtual element scheme, which eventually will imply 
the desired convergence results.

\begin{prop}\label{prop:LPBV}
Let $\bu_h^n=\Bigl(v^{n}_h,w^{n}_h\Bigr)$ be a solution 
of the virtual element scheme \eqref{wf-finite}. 
Then, there exist constants $C>0$, depending
on $\Omega$, $T$, $v_{h}^0$, $w_{h}^0$ and $\Iap$ such that
\begin{equation*}\begin{split}\label{est:L2-norm}
 \|v_{h}\|_{L^\infty(0,T;L^2(\Om))}+\|w_{h}\|_{L^\infty(0,T;L^2(\Om))} \leq C,\\
\Bigl\|\grad v_{h}\Bigr\|_{L^2(\Om_{T})}  \leq C,\\
\|\Pi_k^0 v_{h}\|_{L^4(\Om_{T})} \leq C,
\end{split}
\end{equation*}
where $\Pi_k^0$ has been introduced in \eqref{proyectorg}.
\end{prop}

\begin{proof}
We use (\ref{wf-finite}) with $\varphi_h=v_h^{n}$,
$\phi_h=w^{n}_h$,
and we sum over $n=1,\dots,\kappa$ for all $1<\kappa\leq N$. 
\begin{equation*}\begin{split}
& \sum_{n=1}^\kappa m_h \Bigl(v_{h}^n\!-\!v_{h}^{n-1},v^n_h\Bigl)
+\sum_{n=1}^\kappa m_h \Bigl(w_{h}^n\!-\!w_{h}^{n-1},w^n_h\Bigl)
+\int_0^{\kappa\delt} \D \left( J(v_{h}^n) \right)a_h\Bigl(v_{h}^n,v_h^n\Bigr)\\ & 
\qquad  + \int_0^{\kappa\delt}b_h( v_{h}^n,w_h^n,v_{h}^n) 
= 
\int_0^{\kappa\delt}c_h(v_{h}^n,w_{h}^n,w_{h}^n)+\int_0^{\kappa\delt}\left( {\Iap}_{,h},v_h^n \right)_{0,\Om}.
\end{split}\end{equation*}
Observe that an application of H\"{o}lder and Young inequalities, we get 
\begin{equation*}\begin{split}
&  \sum_{n=1}^{\kappa} m_h \Bigl(v_{h}^n\!-\!v_{h}^{n-1},v^n_h\Bigl)=
\sum_{n=1}^{\kappa}m_h \Bigl(v_{h}^n,v^n_h\Bigl)- \sum_{n=1}^{\kappa}
m_h \Bigl(v_{h}^{n-1},v^n_h\Bigl)\\ 
& \geq  \sum_{n=1}^{\kappa} m_h\Bigl(v_{h}^n,v^n_h\Bigl) -  \sum_{n=1}^{\kappa} \Biggl(m_h\Bigl(v_{h}^n,v^n_h\Bigl)\Biggl)^{1/2}\Biggl(m_h\Bigl(v_{h}^{n-1},v^{n-1}_h\Bigl)\Biggl)^{1/2}\\ 
& \geq  \sum_{n=1}^{\kappa}m_h\Bigl(v_{h}^n,v^n_h\Bigl) - \frac{1}{2} \sum_{n=1}^{\kappa} m_h\Bigl(v_{h}^n,v^n_h\Bigl)-
 \sum_{n=1}^{\kappa}\frac{1}{2}m_h\Bigl(v_{h}^{n-1},v^{n-1}_h\Bigl)
\\ 
&= \sum_{n=1}^{\kappa} \Biggl(\frac{1}{2}m_h\Bigl(v_{h}^n,v^n_h\Bigl)-\frac{1}{2}m_h\Bigl(v_{h}^{n-1},v^{n-1}_h\Bigl)\Biggl)
\\ 
&= \frac{1}{2}m_h\Bigl(v_{h}^\kappa,v_h^\kappa\Bigl)-\frac{1}{2}m_h\Bigl(v_{h}^{0},v^{0}_h\Bigl).
\end{split}\end{equation*}
Using the last inequality, the definition of the forms
$b_h$, $c_h$, the assumption \eqref{Cond_dcap3} and \eqref{stability} we get
\begin{equation*}\begin{split}
& \frac 12 \beta'(v^{\kappa}_h,v^{\kappa}_h)_{0,\Om}
+\frac 12 \beta' (w^{\kappa}_h,w^{\kappa}_h)_{0,\Om}
+d_1 \alpha' \int_0^{\kappa\delt} a\Bigl(v_{h}^n,v_h^n\Bigr)
 + \int_0^{\kappa\delt}\left(\sum_{K\in\mathcal{T}_h}\int_{K} \Iona(\Pi^0_{K,k}v_{h}^n)\Pi^0_{K,k}v_{h}^n\right) 
\\ &   \qquad \leq \frac 12  \beta'' (v^{0}_h,v^{0}_h)_{0,\Om}
+\frac 12  \beta'' (w^{0}_h,w^{0}_h)_{0,\Om}
+\int_0^{\kappa\delt}\left(\sum_{K\in\mathcal{T}_h}\int_{K} H(\Pi^0_{K,k}v_{h}^n,\Pi^0_{K,k}w_{h}^n)\Pi^0_{K,k}w_{h}^n\right)\\ &   
\qquad \qquad-\int_0^{\kappa\delt}\left(\sum_{K\in\mathcal{T}_h}\int_{K} \Ionb(\Pi^0_{K,k}w_{h}^n)\Pi^0_{K,k}v_{h}^n\right) +\int_0^{\kappa\delt}\left( {\Iap}_{,h},v_h^n \right)_{0,\Om}.
\end{split}\end{equation*}
Now, using the definition of bilinear form $a(\cdot,\cdot)$,
and \eqref{eq:ion1}(a) on the left hand side.
Moreover, we use \eqref{eq:ion1}(b), \eqref{eq:ion1}(c)
and Cauchy-Schwarz inequality and the fact that $\Iap(x,t) \in L^{2}(\Om_{T})$,
on the right hand side, we obtain
\begin{equation*}\begin{split}
& \frac 12 \beta'\Vert v^{\kappa}_h\Vert_{0,\Om}^2
+\frac 12 \beta' \Vert w^{\kappa}_h\Vert_{0,\Om}^2
+d_1 \alpha' \int_0^{\kappa\delt} \vert v_{h}^n\vert_{1,\Om}^2
 + \int_0^{\kappa\delt}\left(\sum_{K\in\mathcal{T}_h}\int_{K} \frac{1}{\alpha_1}\vert\Pi^0_{K,k}v_{h}^n\vert^4\right) 
\\ &   \qquad \leq \frac 12  \beta'' \Vert v^{0}_h\Vert_{0,\Om}^2
+\frac 12  \beta'' \Vert w^{0}_h\Vert_{0,\Om}^2
+\int_0^{\kappa\delt}\left(\sum_{K\in\mathcal{T}_h}\int_{K} \vert\Pi^0_{K,k}v_{h}^n\vert\vert\Pi^0_{K,k}w_{h}^n\vert
+\vert\Pi^0_{K,k}w_{h}^n\vert^2\right)\\
&\qquad+\int_0^{\kappa\delt}\left(\sum_{K\in\mathcal{T}_h}\int_{K}
\vert\Pi^0_{K,k}w_{h}^n\vert\vert\Pi^0_{K,k}v_{h}^n\vert\right) +\int_0^{\kappa\delt}\Vert v_h^n\Vert_{0,\Om}^2 + C.
\end{split}\end{equation*}
An application of the Cauchy-Schwarz and Young inequalities, the continuity of
$\Pi^0_{K,k}$ with respect to $\Vert\cdot\Vert_{0,K}$, yields
\begin{equation}\begin{split}\label{eq:equality-estimates}
& \frac 12 \beta'\Vert v^{\kappa}_h\Vert_{L^2(\Om)}^2
+\frac 12 \beta' \Vert w^{\kappa}_h\Vert_{L^2(\Om)}^2
+d_1 \alpha' \int_0^{\kappa\delt} \vert v_{h}^n\vert_{H^1(\Om)}^2
 + \int_0^{\kappa\delt}\left(\sum_{K\in\mathcal{T}_h}\int_{K} \frac{1}{\alpha_1}\vert\Pi^0_{K,k}v_{h}^n\vert^4\right) 
\\ &   \qquad \leq \frac 12  \beta'' \Vert v^{0}_h\Vert_{L^2(\Om)}^2
+\frac 12  \beta'' \Vert w^{0}_h\Vert_{L^2(\Om)}^2+\int_0^{\kappa\delt}\Vert v_h^n\Vert_{L^2(\Om)}^2
+\int_0^{\kappa\delt}\Vert w_h^n\Vert_{L^2(\Om)}^2+C\\
&\qquad\leq \frac 12  \beta'' \Vert v^{0}_h\Vert_{L^2(\Om)}^2
+\frac 12  \beta'' \Vert w^{0}_h\Vert_{L^2(\Om)}^2+\Vert v_h^n\Vert_{L^2(\Om_T)}^2
+\Vert w_h^n\Vert_{L^2(\Om_T)}^2+C,
\end{split}\end{equation}
thus, for some constants $C_1,C_2,C_3>0$. This implies
\begin{equation}\begin{split}\label{eq:equality-estimates1} 
 \frac 12 \beta'\norm{v^{\kappa}_h}_{L^2(\Om)}+\frac 12 \beta'\norm{w^{\kappa}_h}_{L^2(\Om)}
\leq  C_1\norm{v_{h}}_{L^2(\Om_{T})}^2+ C_2\norm{w_{h}}_{L^2(\Om_{T})}^2+C_3.
\end{split}\end{equation}
Therefore, by the discrete Gronwall inequality, yields from \eqref{eq:equality-estimates1}
\begin{align}
\label{est:Linty-L2-1}
\|v_{h}\|_{L^\infty(0,T;L^2(\Om))}+\|w_{h}\|_{L^\infty(0,T;L^2(\Om))} \leq C_4,
\end{align}
for some constant $C_4>0$.
Finally, using \eqref{est:Linty-L2-1} in \eqref{eq:equality-estimates} and \eqref{eq:ion1}, we get 
\begin{align}
\label{est:Linty-L2}
 \|\Pi_k^0v_{h}\|_{L^4(\Om_{T})}+ \Bigl\|\grad v_{h}\Bigr\|_{L^2(\Om_{T})}  \leq C_5,
\end{align}
for some constant $C_5>0$.
This concludes the proof of Lemma \ref{prop:LPBV}.
\end{proof}

\subsection{Compactness argument and convergence}\label{sec-compactness}
In this section, we will use time continuous approximation of
our discrete solution to obtain compactness in $L^2(\Om_{T})$.
For this, we introduce $\bar v_{h}$ and $\bar w_{h}$ the piecewise affine in
$t$ functions in $W^{1,\infty}([0,T];W_h)$ interpolating the states
$(v^{n}_h)_{n=0,\dots,N}\subset W_h$ and $(w^n_{h})_{n=0,\dots,N} \subset W_h$ at the points $(n\delt
)_{n=0,\dots,N}$. Then, we have
\begin{equation}
	\label{eq:timespace-trans}
	\begin{cases}
& m_h(\pt \bar v_{h}(t),\varphi_h) +\D \left( J( v_h(t))\right) a_h \Bigl( v_h(t),\varphi_h\Bigr)
  +b_h (v_h(t),w_h(t),\varphi_h)=\left( I_{app,h}(t),\varphi_h \right)_{0,\Om},\\
&\displaystyle  m_h(\pt \bar w_h(t),\phi_h)=c_h (v_h(t),w_h(t),\phi_h),
\end{cases}\end{equation}
for all $\varphi_h$ and $\phi_h\in W_h$.

\begin{lem}\label{Space-Time-translate}
There exists a positive constant $C>0$ depending on
$\Omega$, $T$, $v_0$ and $\Iap$ such that 
\begin{equation}\label{Space-translate}
\iint_{\Om_{\br} \times (0,T)}m_h\Bigl(v_{h}(x+\br,t)-v_{h}(x,t),v_{h}(x+\br,t)-v_{h}(x,t)\Bigl)\le C\,|\br|^2 , 
\end{equation}
for all $\br \in \R^2$ with $\Om_{\br}:=\{x\in \Om\,|\, x+\br\in \Om\}$, and
\begin{equation}\label{Time-translate}
\iint_{\Om \times
(0,T-\tau)}m_h\Bigl(v_{h}(x,t+\tau)-v_{h}(x,t),v_{h}(x,t+\tau)-v_{h}(x,t)\Bigl)\dx \dt \le C(\tau+\Delta t),
\end{equation}
for all $\tau\in (0,T)$.
\end{lem}

\begin{proof}
In the first step, we provide the 
{\it proof of estimate \eqref{Space-translate}}.
In this regard, we start with the uniform estimate of space translate of $v_{h}$
from the uniform $L^2(\Om_T)$ estimate of $\Grad v_h$.
Observe that from $L^2(0,T;H^1(\Om))$ estimate of $v_h$, we get easily the bound:
\begin{equation}\label{eq:space-estimate1}
m_h^{\br}\Bigl(v_{h}(x+\br,t)-v_{h}(x,t),v_{h}(x+\br,t)-v_{h}(x,t)\Bigl)\leq 
C \int_0^T\!\!\int_{\Om_{\br}} |v_h(x+\br,\cdot)-v_h(x,\cdot)|^2 \leq \mathcal{C} |\br|^2,
\end{equation}
for some constant $\mathcal{C}>0$, where $m_h^{\br}(\cdot,\cdot)$ is the restriction of the bilinear form
$m_h(\cdot,\cdot)$ on $\Omega_{\br}$. It is clear that the right-hand side in
\eqref{eq:space-estimate1} vanishes as $|\br|\to 0$, uniformly in $h$.

Now, we furnish the proof of estimate \eqref{Time-translate}. 
Observe that for all $t\in[0,T-\tau]$, the function $\varphi^v_h$
such that $\varphi^v_h(x,t)=v_{h}(x,t+\tau)-v_{h}(x,t)$ takes value in $W_h$ for
$(x,t)\in \Om_T$. Therefore, we can use this function as a test function in the
weak formulations (\ref{wf-finite}). Moreover, we previously proved uniform in $h$ bounds
on $v_{h}$ and $\nabla v_h$ in $L^2(\Omega_{T})$ and on $\Pi_k^0v_{h}$
in $L^4(\Omega_{T})$. This implies the analogous bounds for the translates
$\varphi^v_h$ and $\nabla \varphi^v_h$ in $L^2(\Omega\times(0,T-\tau))$
and $\Pi_k^0\varphi^v_h$ in $L^4(\Omega \times (0,T-\tau))$. 
  
We integrate the first approximation equation of (\ref{eq:timespace-trans}) with respect to the time
parameter $s\in [t,t+\tau]$ (with $0<\tau <T$).  In the resulting equations, we take 
the test function as the corresponding translate $\varphi^v_h$. The result is 
\begin{equation*}
	\label{sl1-timetranslate1}
	\begin{split}
		&\int_0^{T-\tau} \int_{\Om}  m_h\Bigl(v_{h}(x,t+\tau)-v_{h}(x,t),v_{h}(x,t+\tau)-v_{h}(x,t)\Bigl) \dx \dt\\
		&\qquad 
		=\int_0^{T-\tau} \int_{\Om} \int_{t}^{t+ \tau}  m_h\Bigl(\partial_s \bar v_{h} (x,s),v_{h}(x,t+\tau)-v_{h}(x,t) \Bigl)\ds \dx \dt\\
		&\qquad =-\int_0^{T-\tau} \int_{\Om} \int_{t}^{t+ \tau}  \D \left( J( v_h(x,s))\right) a_h \Bigl( v_h(x,s),v_{h}(x,t+\tau)-v_{h}(x,t)\Bigr) \ds\dx \dt
		\\
			&\qquad\qquad  -\int_0^{T-\tau} \int_{\Om} \int_{t}^{t+ \tau}  b_h (v_{h}(x,s),w_{h}(x,s),v_{h}(x,t+\tau)-v_{h}(x,t)) \ds\dx \dt \\
		& \qquad \qquad +\int_0^{T-\tau} \int_{\Om} \int_{t}^{t+ \tau}  (I_{app,h},v_{h}(x,t+\tau)-v_{h}(x,t))  \ds\dx \dt
		\\
		&\qquad=I_1+I_2+I_3.
	\end{split}
\end{equation*}
Now, we bound these integrals separately.
For the term $I_1$, we have
\begin{equation*}
	\label{sl1-timetranslate2}
	\begin{split}
		\abs{I_1} &\leq C  \Biggl[\int_0^{T-\tau} 
		\int_{\Om} \Biggl(\int_{t}^{t+ \tau} \abs{\Grad v_{h}(x,s)}^2\ds\Biggl)^{2} \dx \dt\Biggl]^{\frac{1}{2}}
		\times \Biggl[\int_0^{T-\tau} 
		\int_{\Om} \abs{\Grad (v_{h}(x,t+\tau)-v_{h}(x,t)}^2 \dx \dt \Biggl]^{\frac{1}{2}}\\
		& \leq C\, \tau.
	\end{split}
\end{equation*}
for some constant $C>0$. Herein, we used the Fubini theorem (recall that $\displaystyle \int_t^{t+\tau}\,ds=\tau=\int_{s-\tau}^s\,dt$), the H\"older
inequality and the bounds in $L^2$ of  $\Grad v_{h}$.
Keeping in mind the growth bound of the nonlinearity $\Ion$, we apply the H\"older inequality
(with $p=4$, $p'=4/3$ in the ionic current term and with $p=p'=2$ in the other ones) to deduce  
\begin{equation*}
	\label{sl1-timetranslate2-1}
	\begin{split}
		\abs{I_2} &\leq C \Biggl( \Biggl[\int_0^{T-\tau} 
		\int_{\Om} \Biggl(\int_{t}^{t+ \tau} \abs{\Pi_k^0v_{h}(x,s)}^4\ds\Biggl)^{2} \dx \dt\Biggl]^{\frac{3}{4}}
		\times \Biggl[\int_0^{T-\tau} \int_{\Om} \abs{ \Pi_k^0\varphi^v_h(x,t)}^4 \dx \dt \Biggl]^{\frac{1}{4}}\\
		&\qquad +\Biggl[\int_0^{T-\tau} \int_{\Om} \Biggl(\int_{t}^{t+ \tau} \abs{w_{h}(x,s)}^2\ds\Biggl)^{2} \dx \dt\Biggl]^{\frac{1}{2}}
		\times \Biggl[ \int_0^{T-\tau} 
		\int_{\Om} \abs{\varphi^v_h(x,t)}^2 \dx \dt \Biggl]^{\frac{1}{2}}\Biggl)\\
		& \leq C\, \tau, 
	\end{split}
\end{equation*}
for some constant $C>0$, where
we have used that $v_{h},\varphi^v_h$
and $w_{h}$ are uniformly bounded in  $L^2$,
and $\Pi_k^0v_{h},\Pi_k^0\varphi^v_h$ are bounded in $L^4$, and the continuity of
$\Pi^0_{K,k}$ with respect to $\Vert\cdot\Vert_{0,K}$.

Analogously we obtain
$$
\abs{I_3} \leq C\, \tau, 
$$
for some constant $C>0$. Collecting the previous inequalities, we readily deduce
$$ \int_0^{T-\tau}\!\!\int_{\Om} m_h\Bigl(v_{h}(x,t+\tau)-v_{h}(x,t),v_{h}(x,t+\tau)-v_{h}(x,t)\Bigl)
 \leq \mathcal C\,\tau.
$$
Note that, it is easily seen from the definition of ($\bar v_h,\bar w_h)$ and from the discrete weak formulation (\ref{wf-finite}) 
 and estimates in Proposition \ref{prop:LPBV} that
$$
\|\bar v_h - v_h\|^2_{L^2(\Om_{T})}\leq \sum_{n=1}^N \Delta t
\|v_{h}^n-v_{h}^{n-1}\|^2_{L^2(\Om)}\leq \,\mathcal{C}(\Delta t)\to 0 \;\;\text{as $\Delta  t \to 0$}.
$$
This concludes the proof of Lemma \ref{Space-Time-translate}.
\end{proof}

\subsection{Convergence of the virtual element scheme}
\label{sec:conv}

For convergence of our numerical scheme we need the following estimate
\begin{equation}\label{eq-eror-polynome}
\norm{\Pi_k^0 u -u}_{L^2(\Om)} \leq C h^{k+1} \norm{u}_{H^{k+1}(\Om)}\text{ for all $u\in H^{k+1}(\Om)$},
\end{equation}
for some constant $C>0$. This result follows
from standard approximation results (see \cite{BS-2008}).

Note that from Lemma \ref{Space-Time-translate} and the stability condition (\ref{stability}), 
we get
\begin{equation*}\begin{split}
&\iint_{\Om_{\br} \times (0,T)} \abs{v_{h}(x+\br,t)-v_{h}(x,t)}^2\dx \dt 
\le \frac{C}{\beta'}\,|\br|^2 ,
\end{split}\end{equation*}
and
\begin{equation*}
\iint_{\Om \times
(0,T-\tau)}\abs{v_{h}(x,t+\tau)-v_{h}(x,t)}^2dx \dt \le  \frac{C}{\beta'}(\tau+\Delta t).
\end{equation*}
Therefore, the next lemma is a consequence of \eqref{eq-eror-polynome}, Lemma \ref{Space-Time-translate} 
and Kolmogorov's compactness criterion (see, e.g., \cite{Brezis}, Theorem IV.25).
\begin{lem}\label{lem-conv:1}
There exists a subsequence of $\bu_h=(v_{h},w_h)$, not relabeled, such 
that, as $h\to 0$,
\begin{equation}\label{limit-strong}
\begin{split}
&v_{h},\Pi_k^0 v_{h}\to v \text{ strongly in $L^2(\Omega_{T})$ and a.e. in $\Omega_{T}$},\\
&w_{h},\Pi_k^0 w_{h}\to w \text{ weakly in $L^2(\Omega_{T})$ and a.e. in $\Omega_{T}$},\\
&\text{$v_h \rightharpoonup v$ weakly in $L^{2}(0,T; H^{1}(\Omega))$},\\
&\Pi_k^0 v_{h}\rightharpoonup v \text{ weakly in $L^4(\Omega_{T})$}. 
\end{split}
\end{equation}
\end{lem}
Now, we are going to show that the limit functions $\bu:=(v, w)$ constructed
in Lemma~\ref{lem-conv:1} constitute a weak solution of the nonlocal system
defined in \eqref{eq:weak-sol-1}.

For that we let $\varphi \in {\mathcal D}(\Omega \times [ 0, T) )$.
We approximate $\varphi$ by $\varphi_h \in C[0,T;L^2(\Omega)]$ such that
$\varphi_h \vert_{(t^{n-1},t^n)}\in \mathcal{P}_k[t^{n-1},t^n;W_h]$ and $\varphi_h(T)=0$,
where
$\mathcal{P}_k[t^{n-1},t^n;W_h]$ denotes the space of polynomials
of degree $k$ or less having values in $W_h$.

Let $\bu_h:=(v_{h},w_h)$ be the unique solution of the fully discrete
method \eqref{wf-finite}. The proof is based on the convergence to zero
as $h$ goes to zero of each term of the problems.

We start with the convergence of the nonlocal diffusion term.
Observe that
\begin{equation}\label{ineqa}
\begin{split}\displaystyle
\Bigl|\D \Bigl( J( v_h)\Bigl) a_h(v_h,\varphi_h) - \D \Bigl( J( v)\Bigr) a(v,\varphi)\Bigl|  &\leq 
\Bigl|\D \Bigl( J( v)\Bigr)[a_h(v_h,\varphi_h)-a(v,\varphi)]\Bigl|\\
&\qquad +\Bigl|\D \Bigl( J( v_h)\Bigl) - \D \Bigl( J( v)\Bigr) \Bigl||a_h(v_h,\varphi_h)|\\
& := A_1 +A_2.
\end{split}
\end{equation}
For $A_2$, we have
\begin{equation*}
\begin{split}
A_2= \Bigl|\D \Bigl( J( v_h)\Bigl) - \D \Bigl( J( v)\Bigr) \Bigl||a_h(v_h,\varphi_h)|&\leq
C \vert J( v_h)-J( v)\vert |a_h(v_h,\varphi_h)|\\
&\leq C (\Vert v_h-v\Vert_{L^2(\Om)}+\Vert v-\Pi_k^0v\Vert_{L^2(\Om)}) |v_h|_{H^1(\Om)} |\varphi_h|_{H^1(\Om)}\\
&\leq C (\Vert v_h-v\Vert_{L^2(\Om)}+h\vert v\vert_{H^1(\Om)}) |v|_{H^1(\Om)} |\varphi|_{H^1(\Om)}
\end{split}
\end{equation*}
where we have used the assumption \eqref{Cond_dcap3}, the definition of $J( v_h)$ in \eqref{defJ}, then
we add and substract an appropriate polynomial function and finally the continuity of
bilinear form $a_h(\cdot,\cdot)$ in \eqref{continuity}. Thus, using \eqref{limit-strong},
we have that (recall that $\varphi \in {\mathcal D}(\Omega \times [ 0, T) )$)
$$
\lim_{h \to 0} \int_0^T A_2\dt=0.
$$

Now, we bound the term $A_1$ in \eqref{ineqa}.
Using the definition of bilinear form $a_h(\cdot,\cdot)$,
the assumption \eqref{Cond_dcap3}, we have
\begin{equation*}
\begin{split}
A_1=\Bigl|\D \Bigl( J( v)\Bigr)&[a_h(v_h,\varphi_h) - a(v,\varphi)]| 
\leq \Bigl| J( v)\Bigl| \sum_{K \in {\mathcal T}_h}|a^K_h(v_h,\varphi_h) -a^K(v,\varphi)|\\
&\displaystyle \leq \Bigl|J( v)\Bigl| \Biggl[\sum_{K \in {\mathcal T}_h}|a^K(\Pi_{K,k}v_h,\Pi_{K,k}\varphi_h) -a^K(v,\varphi)| +
\sum_{K \in {\mathcal T}_h}|S^K(v_h-\Pi_{K,k}v_h,\varphi_h-\Pi_{K,k}\varphi_h)|\Biggl]\\
&\displaystyle \leq C\Vert v\Vert_{L^2(\Om)}\Biggl[ \sum_{K \in {\mathcal T}_h}|a^K(\Pi_{K,k}v_h-v,\Pi_{K,k}\varphi_h)|
+ \sum_{K \in {\mathcal T}_h} |a^K(v,\Pi_{K,k}\varphi_h-\varphi)|\\
& \qquad \qquad \qquad+\sum_{K \in {\mathcal T}_h}|a^K(v_h-\Pi_{K,k}v_h,\varphi_h-\Pi_{K,k}\varphi_h)|\Biggl],
\end{split}
\end{equation*}
where we have added and substracted $a^K(v,\Pi_{K,k}\varphi_h)$
and used \eqref{Sestab1}. Defining $$\Theta(h)=\displaystyle \sum_{K \in {\mathcal T}_h}|a^K(\Pi_{K,k}v_h-v,\Pi_{K,k}\varphi_h)|$$
Now, using this,
and the Cauchy-Schwarz inequality, we obtain,
\begin{equation*}
\begin{split}
A_1&\displaystyle \leq C\Vert v\Vert_{L^2(\Om)}\Biggl[  \Theta(h)
+ \sum_{K \in {\mathcal T}_h} \vert v\vert_{H^{1}(K)}\vert\Pi_{K,k}\varphi_h-\varphi\vert_{H^{1}(K)}
+\sum_{K \in {\mathcal T}_h}\vert v_h
-\Pi_{K,k}v_h\vert_{H^{1}(K)}\vert\varphi_h-\Pi_{K,k}\varphi_h\vert_{H^{1}(K)}\Biggl],
\end{split}
\end{equation*}
Next, we add and substract an appropriate polynomial $\varphi_{\Pi}$ in the second
term, and we add and substract $\varphi$ in the last term. Thus, we have
\begin{equation*}
\begin{split}
A_1&\displaystyle \leq C\Vert v\Vert_{L^2(\Om)}\Biggl[ \Theta(h)
+ \sum_{K \in {\mathcal T}_h} \vert v\vert_{H^{1}(K)}(\vert\Pi_{K,k}(\varphi_h-\varphi_{\Pi})\vert_{H^{1}(K)}
+\vert\varphi-\varphi_{\Pi}\vert_{H^{1}(K)})\\
& \qquad \qquad \qquad+\sum_{K \in {\mathcal T}_h}\vert v_h
-\Pi_{K,k}v_h\vert_{H^{1}(K)}(\vert\varphi_h-\varphi\vert_{H^{1}(K)}+\vert\varphi-\Pi_{K,k}\varphi_h\vert_{H^{1}(K)})\Biggl]\\
&\displaystyle \leq C\Vert v\Vert_{L^2(\Om)}\Biggl[ \Theta(h)
+ \sum_{K \in {\mathcal T}_h} \vert v\vert_{H^{1}(K)}(\vert\varphi_h-\varphi\vert_{H^{1}(K)}
+\vert\varphi-\varphi_{\Pi}\vert_{H^{1}(K)})\\
& \qquad \qquad \qquad+\sum_{K \in {\mathcal T}_h}\vert v_h\vert_{H^{1}(K)}
(\vert\varphi_h-\varphi\vert_{H^{1}(K)}+\vert\varphi-\varphi_{\Pi}\vert_{H^{1}(K)})\Biggl].
\end{split}
\end{equation*}
Now, using \eqref{limit-strong}, standard approximation results
for polynomials, and the regularity of $\varphi$,
we obtain 
$$
\lim_{h \to 0} \int_0^T A_1 \dt=0.
$$ 
Finally, we get
 \begin{equation*}
\int_0^T \Bigl|\D \Bigl( J( v_h)\Bigl) a_h(v_h,\varphi_h) - \D \Bigl( J( v)\Bigr) a(v,\varphi)\Bigl| \dt \to0
\;\;\text{as}\,\,h\to0.
\end{equation*}
Now, we prove 
\begin{equation}\label{limit-first-term}
\left\vert \int_{0}^{T}m_h({v}_h, \partial_t \varphi_h)
- (  v,\partial_t\varphi)_{0,\Om}\right\vert\to0
\;\;\text{as}\,\,h\to0.
\end{equation}
In fact, using the definition of the bilinear form $m_h(\cdot,\cdot)$, we obtain
\begin{equation*}
\begin{split}
\abs{\int_{0}^{T}m_h({v}_h,\partial_t \varphi_h)
- (v,\partial_t\varphi)_{0,\Om}}\leq &\abs{\sum_{K\in\mathcal{T}_h}(\Pi_K^0{v}_h,\Pi_K^0\partial_t \varphi_h)_{0,K}-({v},\partial_t\varphi)_{0,K}}\\
&\qquad\qquad \qquad+\abs{S_0^K({v}_h-\Pi_K^0{v}_h,\partial_t\varphi_h-\Pi_K^0\partial_t\varphi_h)}\\
 \leq &\abs{\sum_{K\in\mathcal{T}_h}(\Pi_K^0{v}_h-v,\Pi_K^0\partial_t \varphi_h)_{0,K}}+\abs{(v,\Pi_K^0\partial_t \varphi_h-\partial_t\varphi)_{0,K}}\\
&\qquad\qquad \qquad+\abs{S_0^K({v}_h-\Pi_K^0{v}_h,\partial_t\varphi_h-\Pi_K^0\partial_t\varphi_h)}\\
\leq & \norm{v_h-v}_{L^2(\Om)}\norm{\partial_t \varphi}_{L^2(\Om)}+\norm{v}_{L^2(\Om)}\norm{\partial_t \varphi_h-\partial_t \varphi}_{L^2(\Om)}\\
&
+\norm{v_h}_{L^2(\Om)}\Bigl(\norm{\partial_t \varphi_h-\partial_t \varphi_\Pi}_{L^2(\Om)}+\norm{\partial_t \varphi-\partial_t \varphi_\Pi}_{L^2(\Om)}\Bigl).
\end{split}
\end{equation*}
Using this, \eqref{limit-strong}, standard approximation results
for polynomials and the regularity of $\varphi$,
we arrive to \eqref{limit-first-term}. Now, we prove
\begin{equation*}
\int_0^T \left\vert b_h(v_h,w_h,\varphi_h)-(\Ion(v,w),\varphi)_{0,\Om}\right\vert \dt\to0
\;\;\text{as}\,\,h\to0.
\end{equation*}
Using the definition of the form $b_h(\cdot,\cdot,\cdot)$
and the decomposition of the ionic current $\Ion(v,w)$ we have
\begin{equation*}
\begin{split}
&\left\vert b_h(v_h,w_h,\varphi_h)-(\Ion(v,w),\varphi)_{0,\Om}\right\vert
=\left\vert\sum_{K\in\mathcal{T}_h}(\Ion(\Pi_K^0v_h,\Pi_K^0w_h),\Pi_K^0\varphi_h)_{0,K}
-(\Ion(v,w),\varphi)_{0,K}\right\vert\\
&=\left\vert\sum_{K\in\mathcal{T}_h}(\Iona(\Pi_K^0v_h),\Pi_K^0\varphi_h)_{0,K}
+(\Ionb(\Pi_K^0w_h),\Pi_K^0\varphi_h)_{0,K}
-(\Iona(v),\varphi)_{0,K}-(\Ionb(w),\varphi)_{0,K}\right\vert\\
&\le\sum_{K\in\mathcal{T}_h}\vert(\Iona(\Pi_K^0v_h),\Pi_K^0\varphi_h)_{0,K}
-(\Iona(v),\varphi)_{0,K}\vert+\vert(\Ionb(\Pi_K^0w_h),\Pi_K^0\varphi_h)_{0,K}
-(\Ionb(w),\varphi)_{0,K}\vert\\
&=:B_1+B_2.
\end{split}
\end{equation*}
Note that since the function $\Ionb$ is a linear function, we get easily
$$
\int_0^T B_2 \dt \to 0 \text{ as $h$ goes to $0$.}
$$
Now, we turn to the term $B_1$, we have the following estimation
\begin{equation*}
\begin{split}
B_1&\leq \sum_{K\in\mathcal{T}_h}\vert(\Iona(\Pi_K^0v_h),\Pi_K^0\varphi_h)_{0,K}
-(\Iona(\Pi_K^0v_h),\varphi)_{0,K}\vert\\
&\qquad \qquad +\sum_{K\in\mathcal{T}_h}\vert(\Iona(\Pi_K^0v_h),\varphi)_{0,K}
-(\Iona(v),\varphi)_{0,K}\vert\\
&\leq \norm{\varphi_h - \varphi}_{L^\infty(\Omega)} \norm{\Iona(\Pi_K^0v_h)}_{L^1(\Om)}+
Const(v,\Pi_K^0v_h, v_h) \norm{\Pi_K^0v_h - v}_{L^2(\Omega)},
\end{split}
\end{equation*}
where $Const(v,\Pi_K^0v_h, v_h)>0$ is a constant.
This implies that
$$
\int_0^T B_1 \dt\to 0 \text{ as $h$ goes to $0$.}
$$
Similarly, we get 
\begin{equation*}
\int_0^T \left\vert \left( I_{app,h},\varphi_h \right)_{0,\Om}-(\Iap(x,t),\varphi)_{0,\Om}\right\vert\dt \to0
\;\;\text{as}\,\,h\to0.
\end{equation*}

With the above convergences and by density, we are ready to identify the limit $\bu=(v,w)$ as a (weak)
solution of the system \eqref{eq:main2}. Finally, let $\varphi \in L^{2}(0,T;H^{1}(\Omega))\cap L^{4}(\Omega_{T})$
and $\phi \in C([0,T];L^{2}(\Omega))$, then 
by passing to the limit $h \to 0$ in the following weak formulation (with the help of Lemma \ref{lem-conv:1})
\begin{align*}
&-\int_0^T m_h( v_h(t), \partial_{t} \varphi_h)
+ \int_0^T \D \left( J(v_h(t))\right) a_h \Bigl(v_h(t),\varphi_h\Bigr)
  +\int_0^T b_h\Bigl(v_h(t),w_h(t),\varphi_h)=\int_0^T(\Iaph(t),\varphi_h)_{0,\Om}\\
&\displaystyle \int_0^T m_h(\pt w_h(t),\phi_h)=\int_0^T c_h\Bigl(v_h(t),w_h(t),\phi_h),
\end{align*}
we obtain the limit $\bu=(v,w)$ which is a solution of system \eqref{eq:main2} in the
sense of Definition \ref{defi-ws}.

\setcounter{equation}{0}
\section{Error estimates analysis}\label{sec-error}
In this section, an error estimates will be developed
to our model \eqref{eq:main2}. 
For technical reason (because of the nonlinearity of $\Ion$),
we need to relax the assumptions \eqref{eq:ion1}.
For the error estimates analysis, we will use the following
assumption on $\Ion$: we assume that  $\Ion$ is a linear function on $v$ and $w$ satisfying
\begin{equation}\label{eq:ion33}
\begin{split}
&\forall \,s_1,s_2,z_1,z_2\in \R \quad  \abs{\Ion(s_1,z_1)-\Ion(s_2,z_2)}
\leq \alpha_7 \Bigl(\abs{s_1-s_2}+\abs{z_1-z_2}\Bigl),
\end{split}\end{equation}
for some constant $\alpha_7>0$.

First, we introduce the projection
$\mP^h:H^1(\Om) \to W_h$ as the solution of the
following well-posed problem:
 \begin{equation*}\label{eq-error1}
 \begin{cases}
 &\mP^h u \in W_h,\\
 &a_h(\mP^h u,\varphi_h)=a(u,\varphi_h) \text{ for all $\varphi_h \in W_h$}.
 \end{cases}
 \end{equation*}
We have the following lemma, the proof can be found in
\cite[Lemma 3.1]{BBMR2016}.

 \begin{lem}\label{lem-error1}
 Let $u \in H^1(\Om)$. Then, there exist $C,\tilde{C}>0$,
independent of $h$, such that
 \begin{equation*}\label{eq-error2}
 \abs{\mP^h u - u}_{H^1(\Om)} \leq C h^{k} \abs{u}_{H^{k+1}(\Om)},
 \end{equation*}
Moreover, if the domain is convex, then
  \begin{equation*}\label{eq-error3}
 \norm{\mP^h u - u}_{L^2(\Om)} \leq \tilde{C} h^{k+1} \abs{u}_{H^{k+1}(\Om)}.
 \end{equation*}
 \end{lem}

Our main result in this section is the following theorem.
 \begin{thm}\label{thm-error1}
Let $(v,w)$ be the solution of system~\eqref{eq:main2}
and let $(v_h(t),w_h(t))$ be the solution of the problem~\eqref{eq:FEM1}.
Then, for all $t\in (0,T)$, we have
\begin{equation}\label{eq-error4}\begin{split}
&\norm{v_h(\cdot,t) - v(\cdot,t)}_{L^2(\Om)}+\norm{w_h(\cdot,t) - w(\cdot,t)}_{L^2(\Om)}\\
&\qquad  \leq C \Biggl[ \norm{v^0-v_h^{0}}_{L^2(\Om)}+\norm{w^0-w_h^{0}}_{L^2(\Om)}
+h^{k+1} \Biggl( \abs{v^0}_{H^{k+1}(\Om)}+\abs{w^0}_{H^{k+1}(\Om)}\\
&\qquad+\int_0^t\Biggl(\abs{\Iap}_{H^{k+1}(\Om)}+\abs{v}_{H^{k+1}(\Om)}
+\abs{w}_{H^{k+1}(\Om)}+\abs{\pt v}_{H^{k+1}(\Om)}
+\abs{\pt w}_{H^{k+1}(\Om)}\Biggl)\dt\Biggl)\Biggl] \\
&\qquad \qquad  \qquad  \qquad \qquad \qquad  \qquad  \qquad \qquad \qquad  \qquad  \qquad\qquad \qquad  \times
\text{exp}\Biggl(\int_0^t\Bigl(1+\abs{v}_{2}\Bigl)\dt\Biggl),
\end{split}\end{equation}
for some constant $C>0$. Moreover, let $\bu_h^n=\Bigl(v_h^n,w_h^n\Bigr)$
be the virtual element solution generated 
by \eqref{wf-finite}. Then, for $n=1,\dots,N$
\begin{equation}\label{eq-error4-bis}\begin{split}
&\norm{v_h^n - v(\cdot,t_n)}_{L^2(\Om)}+\norm{w_h^n- w(\cdot,t_n)}_{L^2(\Om)}\\
&\qquad  \leq C \Biggl[ \norm{v^0-v_{h}^0}_{L^2(\Om)}+\norm{w^0-w_h^{0}}_{L^2(\Om)}
+\Delta t \int_0^{t_{n}}\Bigl(\abs{\partial^2_{tt}v}+\abs{\partial^2_{tt}w}\Bigl)\dt\\
&\qquad \qquad 
+h^{k+1} \Biggl( \abs{v^0}_{H^{k+1}(\Om)}+\abs{w^0}_{H^{k+1}(\Om)}\\
&\qquad +\int_0^{t_{n}}\Biggl(\abs{\Iap}_{H^{k+1}(\Om)}
+\abs{v}_{H^{k+1}(\Om)}+\abs{w}_{H^{k+1}(\Om)}+\abs{\pt v}_{H^{k+1}(\Om)}
+\abs{\pt w}_{H^{k+1}(\Om)}\Biggl)\dt\Biggl]\\
&\qquad \qquad  \qquad  \qquad \qquad \qquad  \qquad  \qquad \qquad \qquad  \qquad  \qquad\qquad \qquad
\times\text{exp}\Biggl(\int_0^{t_n}\Bigl(1+\abs{v}_{2}\Bigl)\dt\Biggl),
\end{split}\end{equation}
\end{thm} 
 
\begin{proof}
We start with the proof of bound \eqref{eq-error4}.
First note that
$$
U_h(\cdot,t) - U(\cdot,t)=(U_h(\cdot,t) - \mP^h U(\cdot,t))+(\mP^h U(\cdot,t)
- U(\cdot,t))\quad \text{ for $U=v,w$}.
$$
Observe that from Lemma \ref{lem-error1}, we get easily for $U=v,w$
\begin{equation}\label{eq-error5}\begin{split}
\norm{\mP^h U(\cdot,t) - U(\cdot,t)}_{L^2(\Om)}& \leq C h^{k+1} \abs{U}_{H^{k+1}(\Om)}\\
&\leq C h^{k+1} \left( \abs{U_0}_{H^{k+1}(\Om)}
+\int_{0}^t \abs{\pt U(\cdot,s)}_{H^{k+1}(\Om)}\ds\right)\\
&=C h^{k+1} \Bigl( \abs{U_0}_{H^{k+1}(\Om)}+\norm{\pt U}_{L^1(0,t;H^{k+1}(\Om))}\Bigl),
\end{split}\end{equation}
for all $t\in (0,T)$.

Observe that, using the continuous and semidiscrete
problems (cf. \eqref{eq:main2} and \eqref{eq:FEM1}),
the definition of the projector $\mP^h$
and the fact that the derivative with respect to time commutes
with this projector, we obtain 

\begin{equation}\label{eq-error6}\begin{split}
&m_h(\pt (v_h-\mP^h v),\varphi^v_h)+\D \left( J( v_h)\right) a_h((v_h-\mP^h v),\varphi^v_h)\\
&\qquad =\Bigl({\Iap}_{h},\varphi^v_h)_{0,\Om}-b_h(v_h,w_h,\varphi^v_h)
- m_h(\pt \mP^h v,\varphi^v_h)-\D \Bigl( J( v_h)\Bigl) a_h(\mP^h v,\varphi^v_h)\\
&\qquad =\Bigl({\Iap}_{h},\varphi^v_h)_{0,\Om}-b_h(v_h,w_h,\varphi^v_h)
- m_h(\mP^h \pt v,\varphi^v_h)-\D \Bigl( J( v_h)\Bigl) a(v,\varphi^v_h)\\
&\qquad=\Biggl[\Bigl({\Iap}_{h},\varphi^v_h)_{0,\Om}-({\Iap},\varphi^v_h)_{0,\Om}\Biggl]
-\Biggl[b_h(v_h,w_h,\varphi^v_h)-(\Ion(v,w),\varphi^v_h)_{0,\Om}\Biggl]\\
&\qquad \qquad +\Biggl[( \pt v,\varphi^v_h)_{0,\Om}-m_h( \mP^h \pt v,\varphi^v_h)\Biggl]
+ \Biggl[\left( \D \left( J( v)\right)- \D \left( J( v_h)\right)\right)a(v,\varphi^v_h)\Biggl]
\\
&\qquad:=I_1+I_2+I_3+I_4.
\end{split}\end{equation}

Now, we are going to bound each term $I_1,\ldots,I_4$.
Regarding the first term $I_1$, we have
\begin{equation}\label{eq-error7}
I_1=(\Pi_k^0{\Iap}-{\Iap},\varphi^v_h)_{0,\Om}
\leq C h^{k+1} \abs{\Iap}_{H^{k+1}(\Om)}\norm{\varphi^v_h}_{L^{2}(\Om)},
\end{equation}
for some constant $C>0$, where we have used the definition
of ${\Iap}_{h}$.
Next, for $I_2$, using the definition
of the form $b_h(\cdot,\cdot,\cdot)$
and adding and substracting adequate terms,
we have
\begin{equation*}\label{eq-error8}\begin{split}
I_2=&-\Biggl[b_h(v_h,w_h,\varphi^v_h)-(\mP^h \Ion(v,w),\varphi^v_h)_{0,\Om}\Biggl]
-\Biggl[(\mP^h \Ion(v,w),\varphi^v_h)_{0,\Om}-(\Ion(v,w),\varphi^v_h)_{0,\Om}\Biggl]
\\
=&-\Biggl[\sum_{K\in\mathcal{T}_h}(\Ion(\Pi_{K,k}^0v_h,\Pi_{K,k}^0w_h),\Pi_{K.k}^0\varphi^v_h)_{0,K}
-( \Ion(\mP^h v,\mP^h w),\varphi^v_h)_{0,K}\Biggl]\\
&\qquad \qquad \qquad \qquad \qquad \qquad\qquad \qquad \qquad \qquad \qquad \qquad
-\Biggl[(\mP^h \Ion(v,w)-\Ion(v,w),\varphi^v_h)_{0,\Om}\Biggl]\\
=&-\Biggl[\sum_{K\in\mathcal{T}_h}(\Ion(\Pi_{K,k}^0v_h,\Pi_{K,k}^0w_h)
-\Ion(\mP^h v,\mP^h w),\varphi^v_h)_{0,K}\Biggl]
-\Biggl[(\mP^h \Ion(v,w)-\Ion(v,w),\varphi^v_h)_{0,\Om}\Biggl]
\\
\le&C\Biggl[\sum_{K\in\mathcal{T}_h}
(\norm{\Pi_{K,k}^0v_h-\mP^h v}_{L^2(K)}+\norm{\Pi_{K,k}^0w_h
-\mP^h w}_{L^2(K)})\norm{\varphi^v_h}_{L^2(K)}\Biggl]\\
&+C h^{k+1} (\abs{v}_{H^{k+1}(\Om)}+\abs{w}_{H^{k+1}(\Om)})\norm{\varphi^v_h}_{L^2(\Om)}
\\
\leq &C \Biggl(\norm{\Pi_{k}^0v_h-\mP^h v}_{L^2(\Om)}+\norm{\Pi_{k}^0w_h
-\mP^h w}_{L^2(\Om)} \Biggl) \norm{\varphi^v_h}_{L^2(\Om)}\\
&\qquad \qquad \qquad \qquad \qquad \qquad\qquad \qquad \qquad \qquad +
C h^{k+1} (\abs{v}_{H^{k+1}(\Om)}+\abs{w}_{H^{k+1}(\Om)})\norm{\varphi^v_h}_{L^2(\Om)},
\end{split}\end{equation*}
for some constant $C>0$, where we have used
that $\Ion$ is a linear function together with \eqref{eq:ion33},
the properties of projectors $\Pi_k^0$ and $\mP^h$
and finally Lemma~\ref{lem-error1}.

For $I_3$, we use the consistency and stability properties
of the bilinear for $m_h(\cdot,\cdot)$ to get
\begin{equation*}\label{eq-error9}\begin{split}
I_3=&
\sum_{K \in \T_h}\Biggl[( \pt v-\Pi_{K,k}^0 \pt v,\varphi^v_h)_{0,K}
+m^K_h( \Pi_{K,k}^0 \pt v-\mP^h \pt v,\varphi^v_h)\Biggl]
\\
\leq & C \sum_{K \in \T_h} \Biggl[ \norm{\pt v-\Pi_{K,k}^0 \pt v}_{L^2(K)}
+\norm{\Pi_{K,k}^0 \pt v-\mP^h \pt v}_{L^2(K)}\Biggl]
\norm{\varphi^v_h}_{L^2(K)}
\\ \leq & C h^{k+1} \abs{\pt v}_{H^{k+1}(\Om)}\norm{\varphi^v_h}_{L^2(\Om)},
\end{split}\end{equation*}
for some constant $C>0$.
Moreover, by using the assumption on $\D$,
an integration by parts,
the Cauchy-Schwarz inequality
and the continuity of projector $\Pi_k^0$, we get
\begin{equation*}\label{eq-error10}\begin{split}
I_4 \leq & C ( \norm{v-v_h}_{L^2(\Om)}+\norm{v-\Pi_k^0v}_{L^2(\Om)})
\norm{\Delta v}_{L^{2}(\Om)}\norm{\varphi^v_h}_{L^{2}(\Om)},
\end{split}\end{equation*}
for some constant $C>0$.

On the other hand, similarly for $w_h$, we obtain
\begin{equation*}
\begin{split}
&m_h(\pt (w_h-\mP^h w),\varphi^w_h)=\left(c_h(v_h,w_h,\varphi^w_h)
- m_h(\pt \mP^h w,\varphi^w_h)\right)\\
&\qquad =c_h(v_h,w_h,\varphi^w_h)- m_h(\mP^h \pt w,\varphi^v_h)
-(H(v,w),\varphi^w_h)_{0,\Om}+( \pt w,\varphi^w_h)_{0,\Om}\\
& \leq \Biggl[ c_h(v_h,w_h,\varphi^w_h)-(\mP^h H(v,w),\varphi^w_h)_{0,\Om} \Biggl]\\
&\quad  + \Biggl[ (\mP^h H(v,w),\varphi^w_h)_{0,\Om}-(H(v,w),\varphi^w_h)_{0,\Om} \Biggl]
+\Biggl[ (\pt w,\varphi^w_h)-m_h(\pt \mP^h w,\varphi^w_h) \Biggl].
\end{split}
\end{equation*}
Now, assuming that $H$ is a linear function
satisfying \eqref{eq:ion33},
repeating the arguments used to bound $I_2$ and $I_3$,
and using once again the properties of projectors $\Pi_k^0$ and $\mP^h$
and finally Lemma~\ref{lem-error1}, we readily obtain,
\begin{equation}\label{eq-error12}\begin{split}
m_h(\pt (w_h-\mP^h w),\varphi^w_h)\le
&C \Biggl(\norm{\Pi_{k}^0v_h-\mP^h v}_{L^2(\Om)}+\norm{\Pi_{k}^0w_h
-\mP^h w}_{L^2(\Om)} \Biggl) \norm{\varphi^w_h}_{L^2(\Om)}\\
&+C h^{k+1} (\abs{v}_{H^{k+1}(\Om)}+\abs{w}_{H^{k+1}(\Om)}+\abs{\pt w}_{H^{k+1}(\Om)})
\norm{\varphi^w_h}_{L^2(\Om)},
\end{split}\end{equation}
for some constant $C>0$.

Collecting the previous results \eqref{eq-error6}-\eqref{eq-error12},
and using the approximation properties of projectors $\Pi_k^0$ and $\mP^h$,
we get
\begin{equation}\label{eq-error13}\begin{split}
&m_h(\pt (v_h-\mP^h v),\varphi^v_h)+m_h(\pt (w_h-\mP^h w),\varphi^w_h)\\
&\qquad  \leq C h^{k+1} \Biggl[ \abs{\Iap}_{H^{k+1}(\Om)}+\abs{v}_{H^{k+1}(\Om)}
+\abs{w}_{H^{k+1}(\Om)}+\abs{\pt v}_{H^{k+1}(\Om)}+\abs{\pt w}_{H^{k+1}(\Om)}\\
&\qquad \qquad 
+C \Bigl(1+\norm{\Delta v}_{L^{2}(\Om)}\Bigl)\Biggl(\norm{v-v_h}_{L^2(\Om)}+\norm{w-w_h}_{L^2(\Om)} \Biggl)\Biggl](\norm{\varphi_h^v}_{L^2(\Om)}
+\norm{\varphi_h^w}_{L^2(\Om)})
\end{split}\end{equation}
Now, we set $\varphi^v_h:= (v_h-\mP^h v)\in W_h$
and $\varphi^w_h:=(w_h-\mP^h w)\in W_h$ in \eqref{eq-error13}, we deduce
\begin{equation*}\label{eq-error14}\begin{split}
&\frac{1}{2}\frac{d}{dt} \Biggl(m_h(v_h-\mP^h v,v_h-\mP^h v)+m_h(w_h-\mP^h w,w_h-\mP^h w)\Biggl)\\
&\qquad  \leq C h^{k+1} \Biggl[\abs{\Iap}_{H^{k+1}(\Om)}+\abs{v}_{H^{k+1}(\Om)}
+\abs{w}_{H^{k+1}(\Om)}+\abs{\pt v}_{H^{k+1}(\Om)}+\abs{\pt w}_{H^{k+1}(\Om)}\Biggl)\\
&\qquad \qquad 
+C \Bigl(1+\norm{\Delta v}_{L^{2}(\Om)}\Bigl)\Biggl(\norm{v-v_h}_{L^2(\Om)}+\norm{w-w_h}_{L^2(\Om)} \Biggl)\Biggl]\\
&\qquad \qquad \qquad \qquad \qquad \qquad \qquad \qquad \qquad \qquad 
\times \Bigl(\norm{(v_h-\mP^h v)}_{L^2(\Om)}
+\norm{(w_h-\mP^h w)}_{L^2(\Om)}\Bigl).
\end{split}\end{equation*}
Herein, we used the equivalence of the norm $\norm{ \cdot}_h:=m_h(\cdot,\cdot)$ with the $L^2$ norm,
integrating the previous bound on $(0, t)$ and
an application of Gronwall inequality, we get
\begin{equation*}\begin{split}
&(\norm{(v_h-\mP^h v)}_{L^2(\Om)}
+\norm{(w_h-\mP^h w)}_{L^2(\Om)})\\
&\qquad  \leq C \Biggl[ \norm{v_0-v_{0,h}}_0+\norm{w_0-w_{0,h}}_0+h^{k+1} \Biggl( \abs{v_0}_{k+1}+\abs{w_0}_{k+1}\\
&\qquad \qquad+\int_0^t\Biggl({\Iap}_{k+1}+\abs{v}_{k+1}+\abs{w}_{k+1}+\abs{\pt v}_{k+1}+\abs{\pt w}_{k+1}\Biggl)\dt\Biggl] \int_0^t\Bigl(1+\abs{v}_{2}\Bigl)\dt).
\end{split}\end{equation*}
Using this and \eqref{eq-error5}, we get \eqref{eq-error4}.

{\it Proof of \eqref{eq-error4-bis}}:
Similarly to \eqref{eq-error4}, observe that for $n=1,\dots,N$
$$
U^n_h - U(\cdot,t_n)=(U^n_h - \mP^h U(\cdot,t_n))+(\mP^h U(\cdot,t_n) - U(\cdot,t_n))
\quad \text{ for $U=v,w$}.
$$
and from Lemma~\ref{lem-error1}, we get easily for $U=v,w$ and for all $t\in (0,T)$
\begin{equation*}\label{eq-error15}
\norm{\mP^h U(\cdot,t_n) - U(\cdot,t_n)}_{L^2(\Om)} \leq C h^{k+1} \Bigl( \abs{U_0}_{H^{k+1}(\Om)}+\norm{\pt U}_{L^1(0,t;H^{k+1}(\Om))}\Bigl),
\end{equation*}
for some constant $C>0$. Next, we bound the term $(U^n_h - \mP^h U(\cdot,t_n))$ for $U=v, w$. Note that using the continuous and fully discrete problems (cf. \eqref{eq:main2} and \eqref{wf-finite}),
the definition of the projector $\mP^h$, we obtain 
\begin{equation}\label{eq-error16}\begin{split}
&m_h\left(\frac{(v^n_h-\mP^h v(\cdot,t_n))-(v^{n-1}_h-\mP^h v(\cdot,t_{n-1}))}{\Delta t},\varphi^v_h\right)
+\D \left( J( v^n_h)\right) a_h((v^n_h-\mP^h v(\cdot,t_n)),\varphi^v_h)\\
&\qquad =\Bigl(I_{app,h}^n,\varphi^v_h)_{0,\Om}- b_h(v^n_h,w^n_h,\varphi^v_h) - m_h\Biggl(\frac{ \mP^h v(\cdot,t_n)-\mP^h v(\cdot,t_{n-1})}{\Delta t},\varphi^v_h\Biggl)
-\D \left( J( v^n_h)\right) a_h(\mP^h v(\cdot,t_n),\varphi^v_h)\\
&\qquad=(I_{app,h}^n,\varphi^v_h)_{0,\Om}- b_h(v^n_h,w^n_h,\varphi^v_h) - ({\Iap}(\cdot,t_{n}),\varphi^v_h)_{0,\Om}- (\Ion(v(\cdot,t_{n}),w(\cdot,t_{n})),\varphi^v_h)_{0,\Om}
+( \pt v(\cdot,t_{n}),\varphi^v_h)\\
&\qquad \qquad -m_h\Biggl(\frac{ \mP^h v(\cdot,t_{n})-\mP^h v(\cdot,t_{n-1})}{\Delta t},\varphi^v_h\Biggl)
+ \left( \D \left( J( v(\cdot,t_{n}))\right)- \D \left( J( v^n_h)\right)\right)a(v(\cdot,t_{n}),\varphi^v_h)\\
&\qquad=\Biggl[\Bigl(I_{app,h}^n,\varphi^v_h)_{0,\Om}-({\Iap}(\cdot,t_{n}),\varphi^v_h)_{0,\Om}\Biggl]
-\Biggl[b_h(v^n_h,w^n_h,\varphi^v_h)-(\Ion(v(\cdot,t_{n}),w(\cdot,t_{n})),\varphi^v_h)_{0,\Om}\Biggl]\\
&\qquad \qquad \qquad+\Biggl[( \pt v(\cdot,t_{n}),\varphi^v_h)-m_h\Biggl(\frac{ \mP^h v(\cdot,t_{n})-\mP^h v(\cdot,t_{n-1})}{\Delta t},\varphi^v_h\Biggl)\Biggl]
\\
&\qquad \qquad \qquad \qquad \qquad + \Biggl[ \left( \D \left( J( v(\cdot,t_{n}))\right)- \D \left( J( v^n_h)\right)\right)a(v(\cdot,t_{n}),\varphi^v_h)\Biggl]
\\
&\qquad:={\mathcal I}_1+{\mathcal I}_2+{\mathcal I}_3+{\mathcal I}_4.
\end{split}\end{equation}
Now, we will bound the terms ${\mathcal I}_1,\ldots,{\mathcal I}_4$. Note that the first term ${\mathcal I}_1$ can be estimated like \eqref{eq-error7}
\begin{equation*}\label{eq-error17}
{\mathcal I}_1\leq C h^{k+1} \abs{\Iap(\cdot,t_{n})}_{H^{k+1}(\Om)}\norm{\varphi^v_h}_{L^2(\Om)},
\end{equation*}
for some constant $C>0$.
Next, for ${\mathcal I}_2$, using the definition of the form $b_h(\cdot,\cdot,\cdot)$,
adding and substracting adequate terms, we have
\begin{equation*}\label{eq-error18}\begin{split}
{\mathcal I}_2=&-\Biggl[b_h(v_h^n w_h^n,\varphi^v_h)-(\mP^h \Ion(v(\cdot,t_n),w(\cdot,t_n)),\varphi^v_h)_{0,\Om}\Biggl]\\
&\qquad \qquad -\Biggl[(\mP^h \Ion(v(\cdot,t_n),w(\cdot,t_n)),\varphi^v_h)_{0,\Om}-(\Ion(v(\cdot,t_n),w(\cdot,t_n)),\varphi^v_h)_{0,\Om}\Biggl]
\\
=&-\Biggl[\sum_{K\in\mathcal{T}_h}(\Ion(\Pi_{K,k}^0v_h^n,\Pi_{K,k}^0w_h^n),\Pi_{K,k}^0\varphi^v_h)_{0,K}
-( \Ion(\mP^h v(\cdot,t_n),\mP^h w(\cdot,t_n)),\varphi^v_h)_{0,K}\Biggl]\\
&\qquad \qquad -\Biggl[(\mP^h \Ion(v(\cdot,t_n),w(\cdot,t_n))-\Ion(v(\cdot,t_n),w(\cdot,t_n)),\varphi^v_h)_{0,\Om}\Biggl]\\
=&-\Biggl[\sum_{K\in\mathcal{T}_h}(\Ion(\Pi_{K,k}^0v_h^n,\Pi_{K,k}^0w_h^n)
-\Ion(\mP^h v(\cdot,t_n),\mP^h w(\cdot,t_n)),\varphi^v_h)_{0,K}\Biggl]\\
&\qquad \qquad -\Biggl[(\mP^h \Ion(v(\cdot,t_n),w(\cdot,t_n))-\Ion(v(\cdot,t_n),w(\cdot,t_n)),\varphi^v_h)_{0,\Om}\Biggl]
\\
\le&C\Biggl[\sum_{K\in\mathcal{T}_h}
(\norm{\Pi_{K,k}^0v_h^n-\mP^h v(\cdot,t_n)}_{L^2(K)}+\norm{\Pi_{K,k}^0w_h^n
-\mP^h w(\cdot,t_n)}_{L^2(K)})\norm{\varphi^v_h}_{L^2(K)}\Biggl]\\
&\qquad \qquad +C h^{k+1} (\abs{v(\cdot,t_n)}_{H^{k+1}(\Om)}+\abs{w(\cdot,t_n)}_{H^{k+1}(\Om)})\norm{\varphi^v_h}_{L^2(\Om)}
\\
\leq &C \Biggl(\norm{\Pi_{k}^0v_h^n-\mP^h v(\cdot,t_n)}_{L^2(\Om)}+\norm{\Pi_{k}^0w_h^n
-\mP^h w(\cdot,t_n)}_{L^2(\Om)} \Biggl) \norm{\varphi^v_h}_{L^2(\Om)}\\
&\qquad \qquad +C h^{k+1} (\abs{v(\cdot,t_n)}_{H^{k+1}(\Om)}+\abs{w(\cdot,t_n)}_{H^{k+1}(\Om)})\norm{\varphi^v_h}_{L^2(\Om)},
\end{split}\end{equation*}
for some constant $C>0$, where we have used
that $\Ion$ is a linear function together with \eqref{eq:ion33},
the properties of projectors $\Pi_k^0$ and $\mP^h$
and finally Lemma~\ref{lem-error1}.

Regarding ${\mathcal I}_3$, we use the consistency and stability properties of the bilinear form $m_h$ to get
\begin{equation*}\label{eq-error19}\begin{split}
{\mathcal I}_3=& \sum_{K \in \T_h}\Biggl[( \pt v(\cdot,t_n),\varphi^v_h)_{0,K}-m^K_h\Biggl(\frac{ \mP^h v(\cdot,t_n)-\mP^h v(\cdot,t_{n-1})}{\Delta t},\varphi^v_h\Biggl)\Biggl]\\
=& \sum_{K \in \T_h}\Biggl[\Biggl( \pt v(\cdot,t_n)-\frac{ v(\cdot,t_n)- v(\cdot,t_{n-1})}{\Delta t},\varphi^v_h\Biggl)_{0,K}\\
&\qquad \qquad+\Biggl( \frac{ v(\cdot,t_n)- v(\cdot,t_{n-1})}{\Delta t}-\frac{ \Pi_{K,k}^0(v(\cdot,t_n)- v(\cdot,t_{n-1}))}{\Delta t},\varphi^v_h\Biggl)_{0,K}\\
&\qquad \qquad
+m^K_h\Biggl( \frac{ \Pi_{K,k}^0(v(\cdot,t_n)- v(\cdot,t_{n-1}))}{\Delta t}-\frac{ \mP^h(v(\cdot,t_n)- v(\cdot,t_{n-1}))}{\Delta t},\varphi^v_h\Biggl)\Biggl]
\\
\leq & \frac{C}{\Delta t} \sum_{K \in \T_h} \Biggl[ \norm{\Delta t \pt v(\cdot,t_n)-(v(\cdot,t_n)- v(\cdot,t_{n-1}))}_{L^2(K)}\\
&\qquad \qquad
+\norm{(v(\cdot,t_n)- v(\cdot,t_{n-1}))-\Pi_{K,k}^0(v(\cdot,t_n)- v(\cdot,t_{n-1}))}_{L^2(K)}
\\
&\qquad \qquad \qquad \qquad+\norm{\Pi_{K,k}^0(v(\cdot,t_n)- v(\cdot,t_{n-1}))-\mP^h(v(\cdot,t_n)- v(\cdot,t_{n-1}))}_{L^2(K)}\Biggl]
\norm{\varphi^v_h}_{L^2(K)}
\\ \leq & \frac{C}{\Delta t}  \Biggl [  \norm{\Delta t \pt v(\cdot,t_n)-(v(\cdot,t_n)- v(\cdot,t_{n-1}))}_{L^2(\Om)}+h^{k+1} \abs{v(\cdot,t_n)- v(\cdot,t_{n-1})}_{H^{k+1}(\Omega)}\Biggl]
\norm{\varphi^v_h}_{L^2(\Om)}\\ 
\leq & \frac{C}{\Delta t} \Biggl [ \Delta t \int_{t_{n-1}}^{t_n} \norm{\partial^2_{tt} v(\cdot,s)}_{L^2(\Om)}\ds
+h^{k+1} \int_{t_{n-1}}^{t_n}  \abs{v_t(\cdot,s)}_{H^{k+1}(\Om)}\ds\Biggl]
\norm{\varphi^v_h}_{L^2(\Om)}\,
\end{split}\end{equation*}
for some constant $C>0$, where we have used Cauchy-Schwarz inequality,
and the approximation properties of $\Pi_k^0$ and $\mP^h$, finally Lemma~\ref{lem-error1}.
Moreover, for ${\mathcal I}_4$ by using an integration by parts, the assumption on $\D$, Cauchy-Schwarz inequality, and the continuity of projector $\Pi_k^0,$  we get
\begin{equation*}\label{eq-error20-bis}\begin{split}
{\mathcal I}_4\leq & C \left( \norm{v(\cdot,t_n)-v^n_h}_{L^2(\Om)}+ \norm{v(\cdot,t_n)- \Pi^0_k v(\cdot,t_n)}_{L^2(\Om)}\right) \norm{\Delta v(\cdot,t_n)}_{L^2(\Om)}\norm{\varphi^v_h}_{L^2(\Om)},
\end{split}\end{equation*}
for some constant $C>0$.
On the other hand, similarly for $w_h$, we obtain
\begin{equation*}
\begin{split}
& m_h\left(\frac{(w^n_h(\cdot)-\mP^h w(\cdot,t_n))-(w^{n-1}_h(\cdot)-\mP^h w(\cdot,t_{n-1}))}{\Delta t},\varphi^w_h\right) =c_h(v_h^n,w_h^n,\varphi^w_h)\\
&\qquad \qquad\qquad \qquad\qquad \qquad- m_h\left(\frac{ \mP^h w(\cdot,t_n)-\mP^h w(\cdot,t_{n-1})}{\Delta t},\varphi^w_h\right)\\
&\qquad \qquad  =c_h(v_h^n,w_h^n,\varphi^w_h)-m_h\left(\frac{ \mP^h w(\cdot,t_n)-\mP^h w(\cdot,t_{n-1})}{\Delta t},\varphi^w_h\right)\\
&\qquad \qquad \qquad\qquad \qquad  -(H(v(\cdot,t_n),w(\cdot,t_n)),\varphi^w_h)_{0,\Om}+( \pt w(\cdot,t_n),\varphi^w_h)_{0,\Om}\\
& \qquad \qquad \leq \Biggl[ c_h(v_h^n,w_h^n,\varphi^w_h)-(\mP^h H(v(\cdot,t_n),w(\cdot,t_n)),\varphi^w_h)_{0,\Om} \Biggl]\\
&\qquad \qquad \qquad\qquad \qquad + \Biggl[ (\mP^h H(v(\cdot,t_n),w(\cdot,t_n)),\varphi^w_h)_{0,\Om}-(H(v(\cdot,t_n),w(\cdot,t_n)),\varphi^w_h)_{0,\Om} \Biggl]\\
&\qquad \qquad \qquad\qquad \qquad+\Biggl[ ( \pt w(\cdot,t_n),\varphi^w_h)_{0,\Om}-m_h\left(\frac{ \mP^h w(\cdot,t_n)-\mP^h w(\cdot,t_{n-1})}{\Delta t},\varphi^w_h\right) \Biggl].
\end{split}
\end{equation*}
Now, assuming that $H$ is a linear function
satisfying \eqref{eq:ion33},
repeating the arguments used to bound ${\mathcal I}_2$ and ${\mathcal I}_3$,
and using once again the properties of projectors $\Pi_k^0$ and $\mP^h$
and finally Lemma~\ref{lem-error1}, we readily obtain,
\begin{equation}\label{eq-error12-bis}\begin{split}
&m_h\left(\frac{(w^n_h(\cdot)-\mP^h w(\cdot,t_n))-(w^{n-1}_h(\cdot)-\mP^h w(\cdot,t_{n-1}))}{\Delta t},\varphi^w_h\right)\le
C \left(\norm{\Pi_{k}^0v_h^n-\mP^h v(\cdot,t_n)}_{L^2(\Om)}\right.\\
&\qquad \qquad+\norm{\Pi_{k}^0w_h^n
-\mP^h w(\cdot,t_n)}_{L^2(\Om)}\\
&\qquad \qquad+C h^{k+1} \left.(\abs{v(\cdot,t_n)}_{H^{k+1}(\Om)}+\abs{w(\cdot,t_n)}_{H^{k+1}(\Om)})\right)\norm{\varphi^w_h}_{L^2(\Om)}\\
&\qquad \qquad+\frac{C}{\Delta t} \Biggl [ \Delta t \int_{t_{n-1}}^{t_n} \norm{\partial^2_{tt} w(\cdot,s)}_{L^2(\Om)}\ds
+h^{k+1} \int_{t_{n-1}}^{t_n}  \abs{w_t(\cdot,s)}_{H^{k+1}(\Om)}\ds\Biggl]
\norm{\varphi^v_h}_{L^2(\Om)},
\end{split}\end{equation}
for some constant $C>0$.

Collecting the previous results \eqref{eq-error16} - \eqref{eq-error12-bis}, and using the approximation properties of projectors $\Pi_k^0$ and $\mP^h$, after substituting 
$\varphi^v_h= v_h-\mP^h v$ and $\varphi^w_h=w_h-\mP^h w$ in \eqref{eq-error16} and \eqref{eq-error12-bis}, respectively, we deduce
\begin{equation*}\label{eq-error14-bis-bis}\begin{split}
& \Biggl(m_h\Bigl(v^n_h-\mP^h v(\cdot,t_n),v^n_h-\mP^h v(\cdot,t_n)\Bigl)+m_h\Bigl(w^n_h-\mP^h w(\cdot,t_n),w^n_h-\mP^h w(\cdot,t_n)\Bigl)\Biggl)\\
&\quad  \leq \Biggl(m_h(v^{n-1}_h-\mP^h v(\cdot,t_{n-1}),v^n_h-\mP^h v(\cdot,t_n))+m_h(w^{n-1}_h-\mP^h w(\cdot,t_{n-1}),w^n_h-\mP^h w(\cdot,t_n))\Biggl)\\
&\quad \quad  +C \Delta t \Biggl[(1+\abs{v(\cdot,t_n)}_{2})\norm{\Pi_k^0v^n_h-\mP^h v(\cdot,t_n)}_{L^2(\Om)}+\norm{\Pi_k^0w^n_h-\mP^h w(\cdot,t_n)}_{L^2(\Om)} \\
&\quad \quad \quad+
h^{k+1} \Bigl(\abs{v(\cdot,t_n)}_{H^{k+1}(\Om)}+\abs{w(\cdot,t_n)}_{H^{k+1}(\Om)}\Bigl)
+  \int_{t_{n-1}}^{t_n} \Bigl(\norm{\partial^2_{tt} v(\cdot,s)}_{L^2(\Om)}+\norm{\partial^2_{tt} w(\cdot,s)}_{L^2(\Om)}\Bigl)\ds\\
&\qquad \qquad 
+h^{k+1} \int_{t_{n-1}}^{t_n}  \Bigl(\abs{v_t(\cdot,s)}_{H^{k+1}(\Om)}+\abs{w_t(\cdot,s)}_{H^{k+1}(\Om)}\Bigl)\ds+ h^{k+1} \abs{\Iap(\cdot,t_n)}_{H^{k+1}(\Om)}\Biggl]\times \\
&\qquad \qquad \quad
 \Biggl(\norm{v^n_h-\mP^h v(\cdot,t_n))}+\norm{w^n_h-\mP^h w(\cdot,t_n)}\Biggl).
\end{split}\end{equation*}

This implies
\begin{equation}\label{eq-error-est-final1}\begin{split}
&\Biggl(\norm{v^n_h-\mP^h v(\cdot,t_n)}_h+\norm{w^n_h-\mP^h w(\cdot,t_n)}_h\Biggl)\\
&\quad  \leq \Biggl(\norm{v^{n-1}_h-\mP^h v(\cdot,t_{n-1})}_h+\norm{w^{n-1}_h-\mP^h w(\cdot,t_{n-1})}_h\Biggl)\\
&\quad \quad  +C \Delta t \Biggl[(1+\abs{v(\cdot,t_n)}_{2})\norm{\Pi_{k}^0 v^n_h-\mP^h v(\cdot,t_n)}_{L^2(\Om)}+\norm{\Pi_{k}^0 w^n_h-\mP^h w(\cdot,t_n)}_{L^2(\Om)} \\
&\quad \quad \quad+
h^{k+1} \Bigl(\abs{v(\cdot,t_n)}_{H^{k+1}(\Om)}+\abs{w(\cdot,t_n)}_{H^{k+1}(\Om)}\Bigl)
+  \int_{t_{n-1}}^{t_n} \Bigl(\norm{\partial^2_{tt} v(\cdot,s)}_{L^2(\Om)}+\norm{\partial^2_{tt} w(\cdot,s)}_{L^2(\Om)}\Bigl)\ds\\
&\qquad \qquad 
+h^{k+1} \int_{t_{n-1}}^{t_n}  \Bigl(\abs{v_t(\cdot,s)}_{H^{k+1}(\Om)}+\abs{w_t(\cdot,s)}_{H^{k+1}(\Om)}\Bigl)\ds+ h^{k+1} \abs{\Iap(\cdot,t_n)}_{H^{k+1}(\Om)}\Biggl]
\end{split}\end{equation}
\begin{equation*}\label{eq-error-est-final2}\begin{split}
&\quad  \leq \Biggl(\norm{v^{0}_h-\mP^h v(\cdot,0)}_h+\norm{w^{0}_h-\mP^h w(\cdot,0)}_h\Biggl)\\
&\qquad \qquad  +
C \Delta t \sum_{l=1}^n\Biggl[(1+\abs{v(\cdot,t_l)}_{2})\Bigl(\norm{v^l_h-\mP^h v(\cdot,t_l)}_h+\norm{w^l_h-\mP^h w(\cdot,t_l)}_h \Bigl)\\
&\qquad \qquad
+  \int_{t_{l-1}}^{t_l} \Bigl(\norm{\partial^2_{tt} v(\cdot,s)}_{L^2(\Om)}+\norm{\partial^2_{tt} w(\cdot,s)}_{L^2(\Om)}\Bigl)\ds
+\norm{v_0-v_{0,h}}_{L^2(\Om)}+\norm{w_0-w_{0,h}}_{L^2(\Om)}\\
&\qquad \qquad 
+h^{k+1} \Biggl( \abs{v_0}_{H^{k+1}(\Om)}+\abs{w_0}_{H^{k+1}(\Om)}+\abs{v(\cdot,t_l)}_{H^{k+1}(\Om)}+\abs{w(\cdot,t_l)}_{H^{k+1}(\Om)}\\
&\qquad \qquad \qquad +\int_{t_{l-1}}^{t_l}  \Bigl(\abs{v(\cdot,s)}_{H^{k+1}(\Om)}+\abs{w(\cdot,s)}_{H^{k+1}(\Om)}\Bigl)\ds
+ \abs{\Iap(\cdot,t_l)}_{H^{k+1}(\Om)}\Biggl)
\end{split}\end{equation*}
Finally, we use the equivalence of the norm $\norm{ \cdot}_h:=m_h(\cdot,\cdot)$ with the $L^2$ norm
and an application of discrete Gronwall inequality to \eqref{eq-error-est-final1} to get \eqref{eq-error4-bis}.
This concludes the proof of Theorem \ref{thm-error1}.

\end{proof}

\setcounter{equation}{0}
\section{Numerical results}
\label{sec:numerical-results}

In the present section, we report some numerical examples
of the proposed virtual element method.
With this aim, we have implemented in a MATLAB code the
lowest-order VEM ($k=1$) on arbitrary
polygonal meshes following the ideas proposed in \cite{BBMR2014}.
Moreover, we solve the nonlinear problem derived from \eqref{wf-finite}
by a classical Picard-type iteration.

To complete the choice of the VEM, we have to choose  the bilinear forms
$S^K(\cdot,\cdot)$ and $S_{0}^K(\cdot,\cdot)$ satisfying
\eqref{Sestab1} and \eqref{Sestab2}, respectively. In this respect, we
have proceeded as in \cite[Section 4.6]{Beirao2}: for each polygon $K$ with
vertices $P_1,\dots,P_{N_{K}}$, we have used
\begin{align*}
S^{K}(u,v)&:=\sum_{r=1}^{N_{K}}u(P_r)v(P_r),
\qquad u,v\in W_{1|K},\\
S_{0}^{K}(u,v)&:=h_{K}^{2}\sum_{r=1}^{N_{K}}u(P_r)v(P_r),
\qquad u,v\in W_{1|K}.
\end{align*}
A proof of \eqref{Sestab1}-\eqref{Sestab2} for the above
(standard) choices could be derived following the
arguments in \cite{AABMR13,Beirao2,BLR2017}.
The choices above are standard in the Virtual
Element Literature and correspond to a scaled identity matrix
in the space of the degrees of freedom values.

\begin{center}
\begin{figure}[ht]
\begin{tabular}{ccc}
\includegraphics[width=5cm, height=5cm]{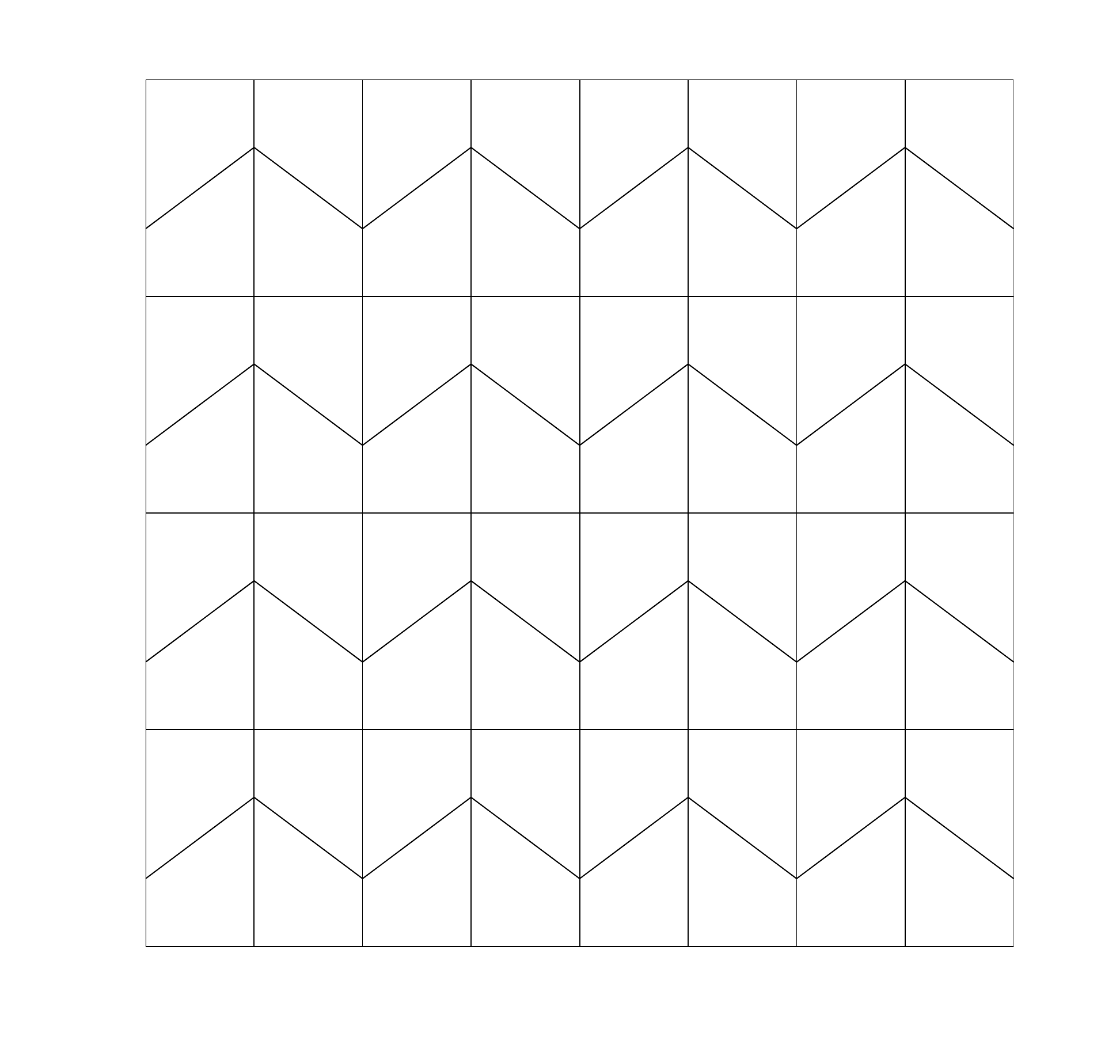} &
\includegraphics[width=5cm, height=5cm]{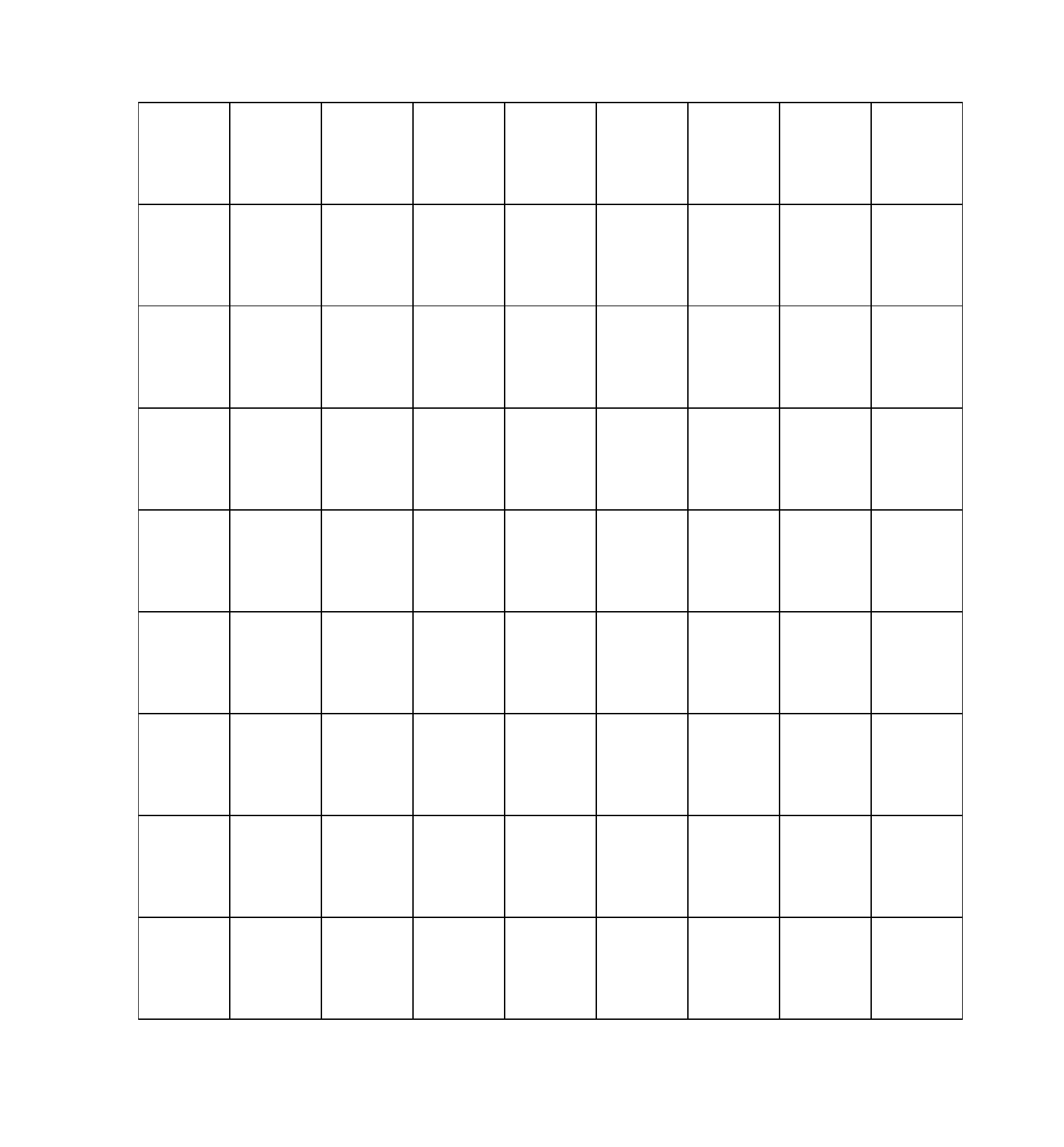}&
\includegraphics[width=5cm, height=5cm]{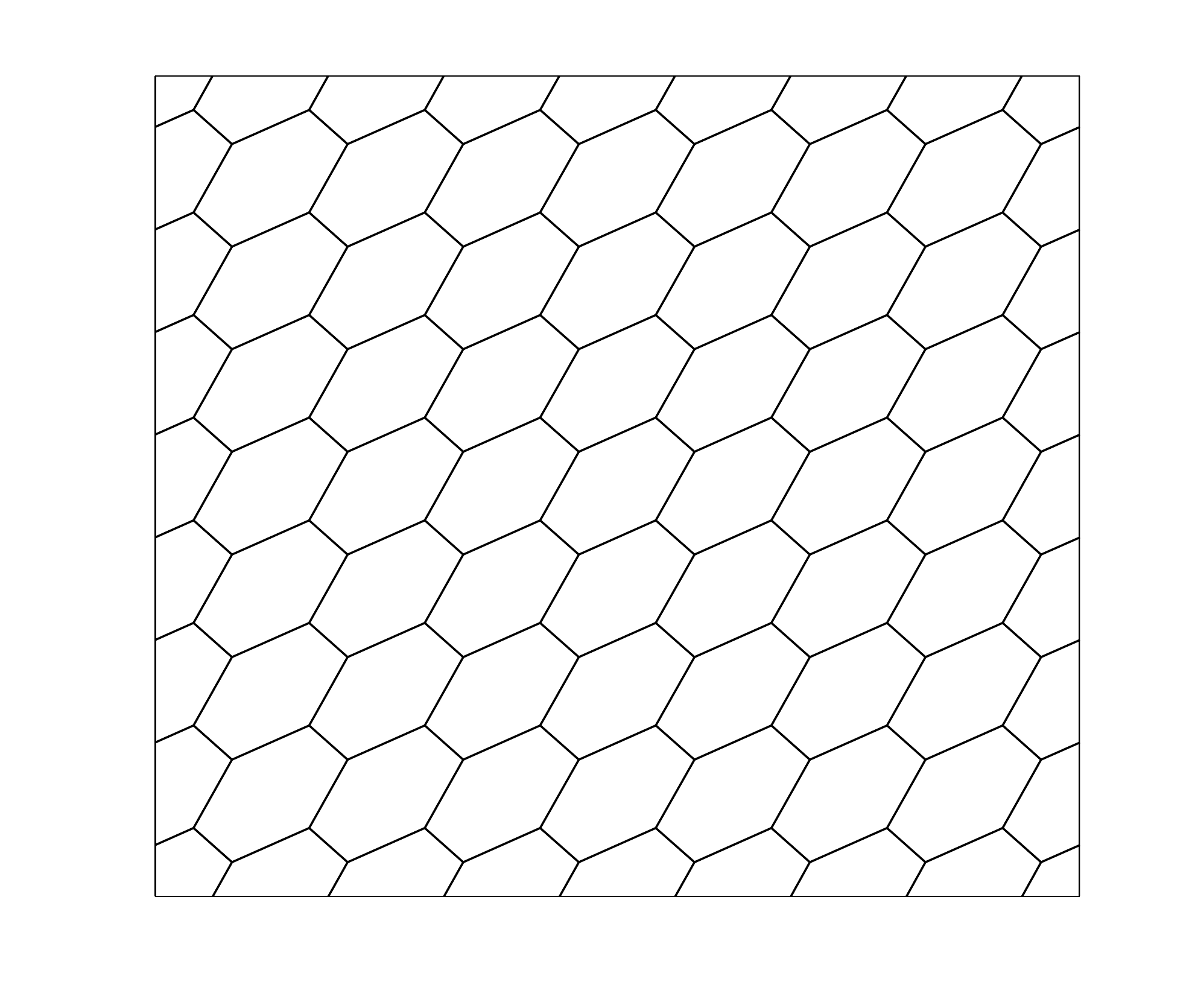}\\
\end{tabular}
\caption{Sample meshes: ${\mathcal T}_h^1$ (left), ${\mathcal T}_h^2$ (center),
${\mathcal T}_h^3$ (right).}
\label{malla}
\end{figure}
\end{center}

In all the numerical examples we have considered $H(v,w)$
and $\Ion(v,w)$ as in \eqref{eq:H} and \eqref{eq:ionplus}, respectively.
Moreover, we have tested the method by using
different families of meshes
(see Figure~\ref{malla}).

\subsection{Example 1}

The aim of this numerical example is to test
the convergence properties of the proposed VEM.
With this objective, we introduce the following
discrete relative $L^{2}$ norm of the difference
between the interpolant $w_I$ of a reference solution
obtain on an extremely fine mesh
and the numerical solution $w_h$ at the final time
$T$, that is,
$$E_{h,\Delta t}^2:=\frac{m_h(w_I(\cdot,T)-w_h(\cdot,T),w_I(\cdot,T)-w_h(\cdot,T))}
{m_h(w_I(\cdot,T),w_I(\cdot,T))}.$$

For this example, the domain will be $\Omega=(0,1)^2$
and the time interval $[0,1]$,
we will take the model constants as follow $a=0.2232, b=0.9, \lambda=-1, \theta=0.004$.
We also take $\Iap=0$ and $D(x)=0.01x$.
Moreover, we consider the following initial data:
\begin{equation*}
v^0(x,y)=\Biggl(1+0.5\cos(4\pi x)\cos(4\pi y) \Biggl),
\qquad w^0(x,y)=\Biggl(1+0.5\cos(8\pi x)\cos(8\pi y) \Biggl).
\end{equation*}

Due to the lack of exact solution for this example, we compute errors
using a numerical solution on an extremely fine mesh $(h=1/512)$ and
time step $(\Delta t=1/200)$ as reference $v_{ref}, \, w_{ref}$.

We report in Table~\ref{eru2} the relative errors
in the discrete $L^{2}$-norm of the variable $v$,
for the family of meshes ${\mathcal T}_h^2$
and different refinement levels and time steps.

\begin{center}
\begin{table}[ht]
\begin{tabular}{|l|c|c|c|c|}
\hline 
$h\backslash\Delta t$&  $\Delta t=1/10$ &$\Delta t=1/20$ & $\Delta t=1/40$ & $\Delta t=1/80$\\
\hline  
 $1/8$   &$0.039090164250364$& $0.022444307847497$ &$0.017875789409797$&$0.016570898457504$\\ 
 $1/16$  &$0.032397693528292$& $0.013142430630831$ &$0.007095750525537$&$0.005512035102741$\\
 $1/32$ &$0.031718009281784$& $0.011763920419019$ &$0.004646804404667$&$0.002076778899273$\\
 $1/64$  &$0.031626299496412$& $0.011596060410183$ &$0.004354704245183$&$0.001528339183782$\\
 \hline 
\end{tabular}
\caption{Test 1: Computed error in the discrete $L^2$ norm for $v$.}
\label{eru2}
\end{table}
\end{center}

It can be seen from Table~\ref{eru2} 
that the error in the discrete $L^2$ norm
reduced with a quadratic order with respect to $h$,
which is the expected order of convergence for $k=1$.

\begin{figure}[ht]
\begin{center}
\begin{tabular}{cc}
\includegraphics[width=5.5cm, height=5.5cm]{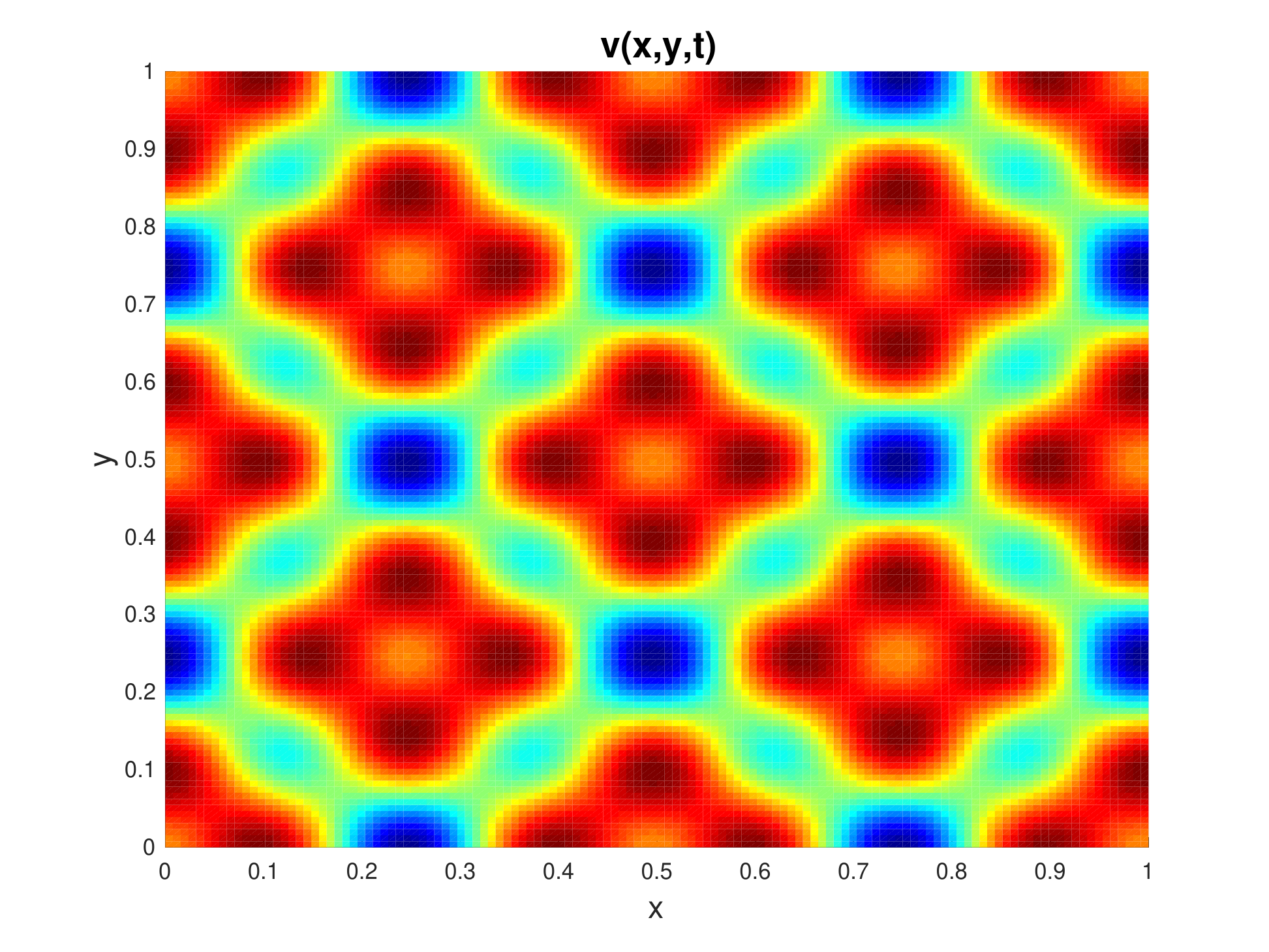} &
\includegraphics[width=5.5cm, height=5.5cm]{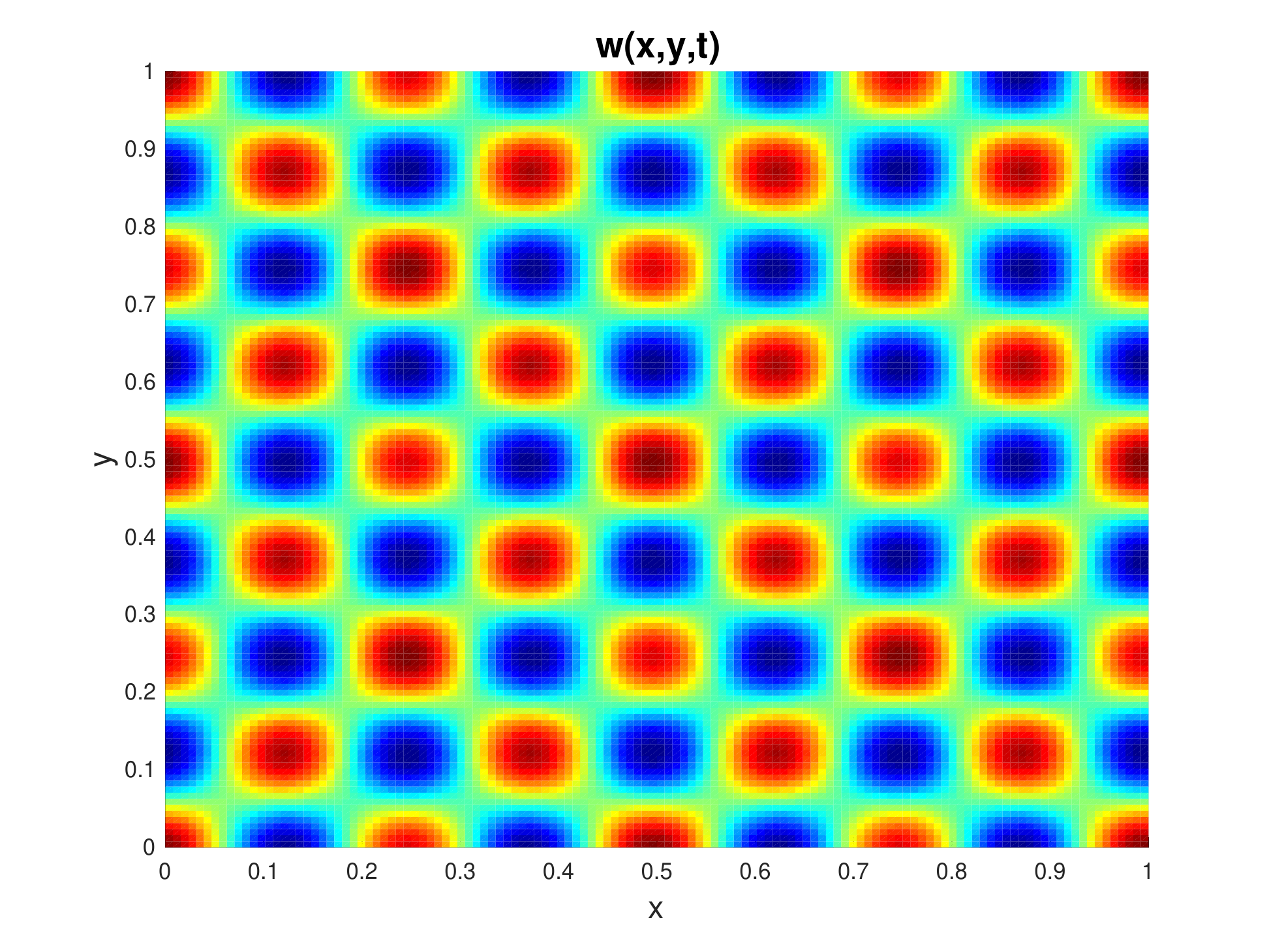} 
\end{tabular}
\end{center}
\caption{Variables $v$ and $w$ for $h=1/64$ and $\Delta t=1/80$.}
\label{fig4}
\end{figure}

We show in Figure~\ref{fig4} the profiles of the
computed quantities.

\subsection{Example 2}

In this test, we consider a benchmark example.
We solve the FitzHugh-Nagumo
equation using meshes $\mathcal T_h^1$ (with $h=1/128$) on the unit square
and time interval $[0,5]$ (with $\Delta t=1/100$)
and with the following 
model constants: $a=0.16875, b=1, \lambda=-100, \theta=0.25$.
Moreover, we consider the following initial data:
\begin{equation*}
v^0(x,y)=\Biggl(1-\frac{1}{1+e^{-50(x^2+y^2)^{1/2}-0.1}}\Biggl), \qquad w^0(x,y)=0.
\end{equation*}
After $4ms$, an instantaneous stimulus is applied in $(x_0,y_0)=(0.5,0.5)$
to the transmembrane potential $v$,
$$I_{app}=\left\{\begin{array}{ll}
 1 \quad mV & \text{if} \qquad (x-x_0)^2+ (y-y_0)^2< 0.04 \qquad cm^2,\\
0 \quad mV & \text{otherwise}.
\end{array}
\right.$$

\begin{figure}[p]
\begin{tabular}{ccc}
\includegraphics[width=5cm, height=5cm]{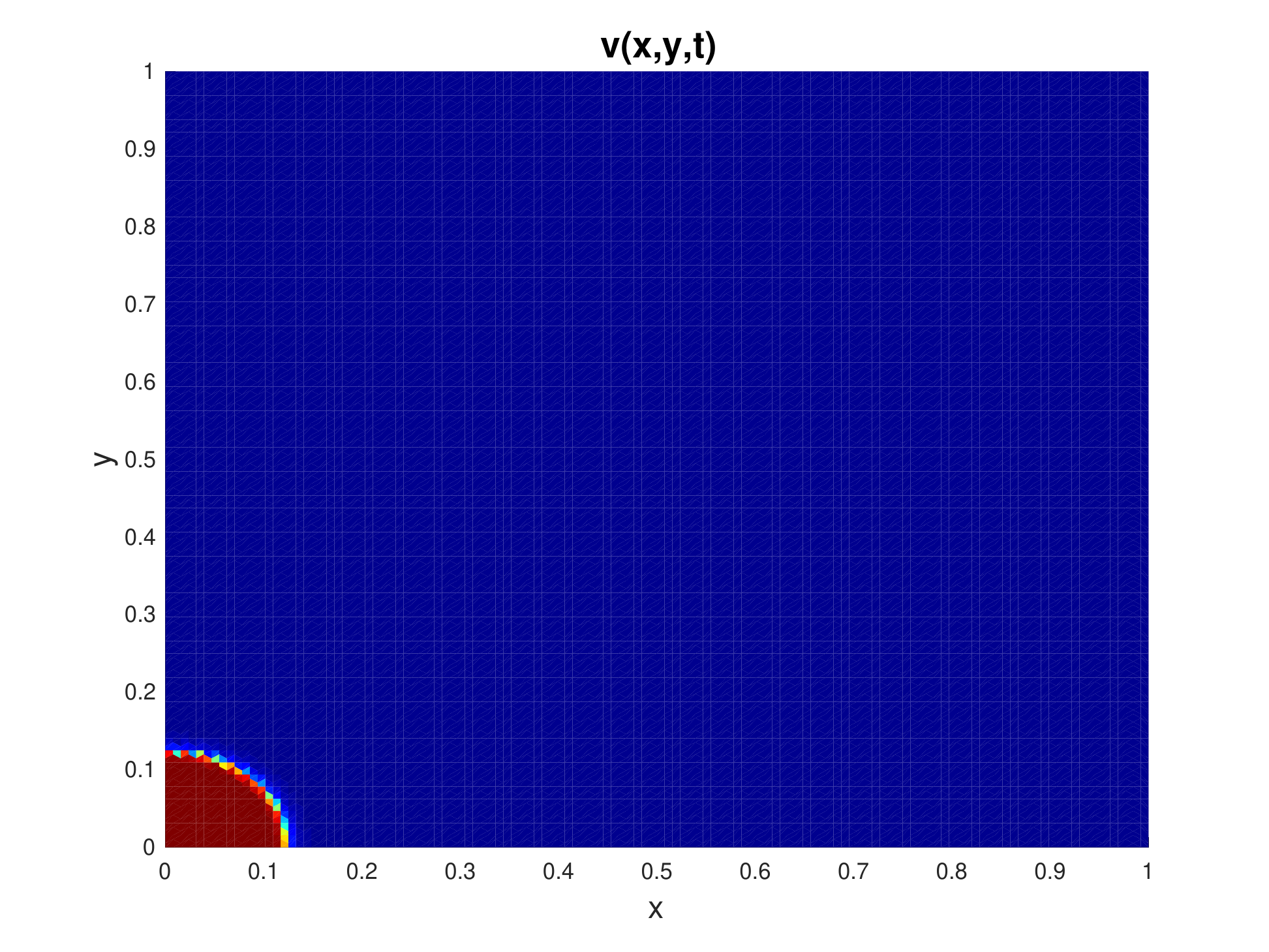} &
\includegraphics[width=5cm, height=5cm]{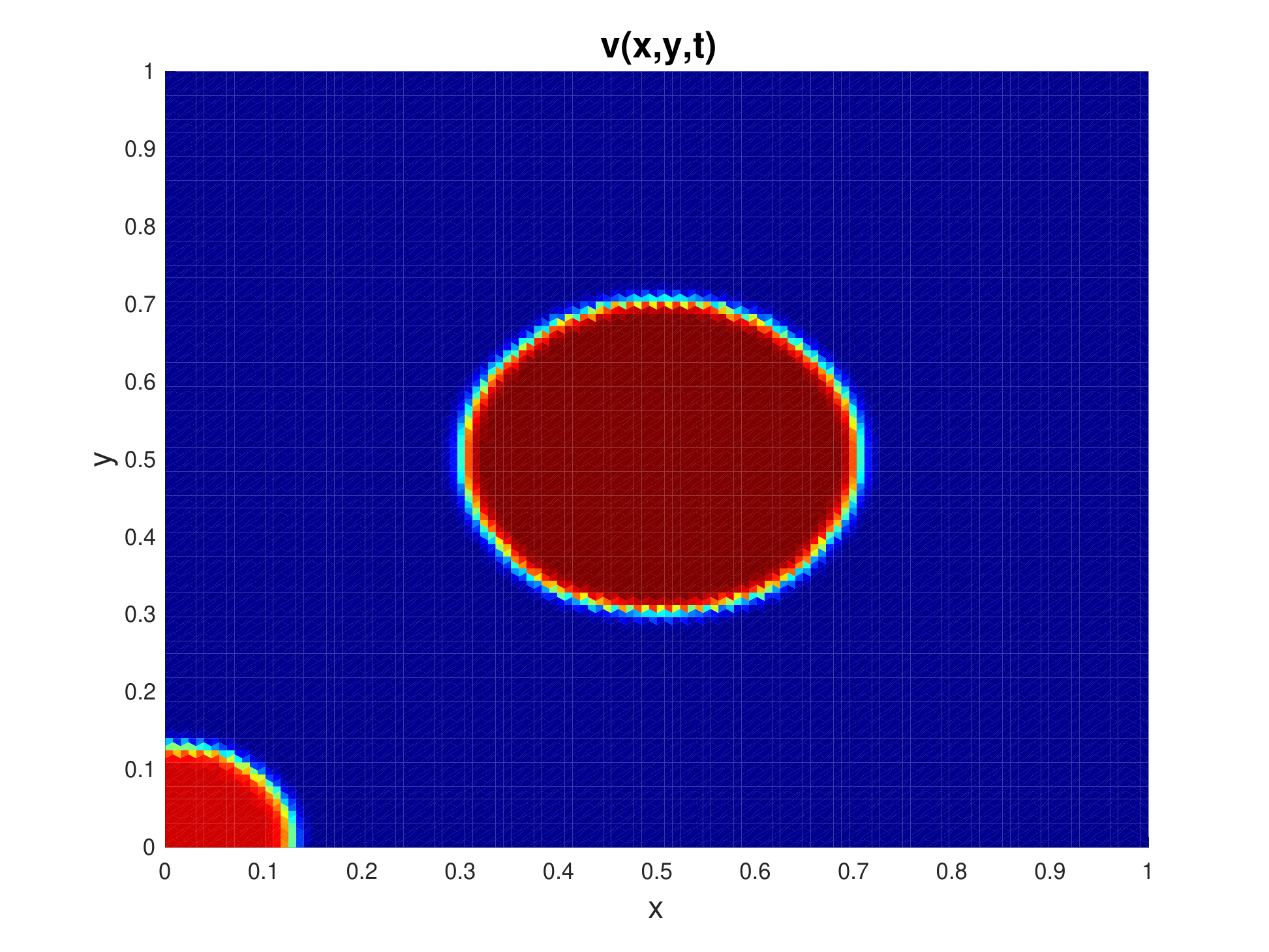} &
\includegraphics[width=5cm, height=5cm]{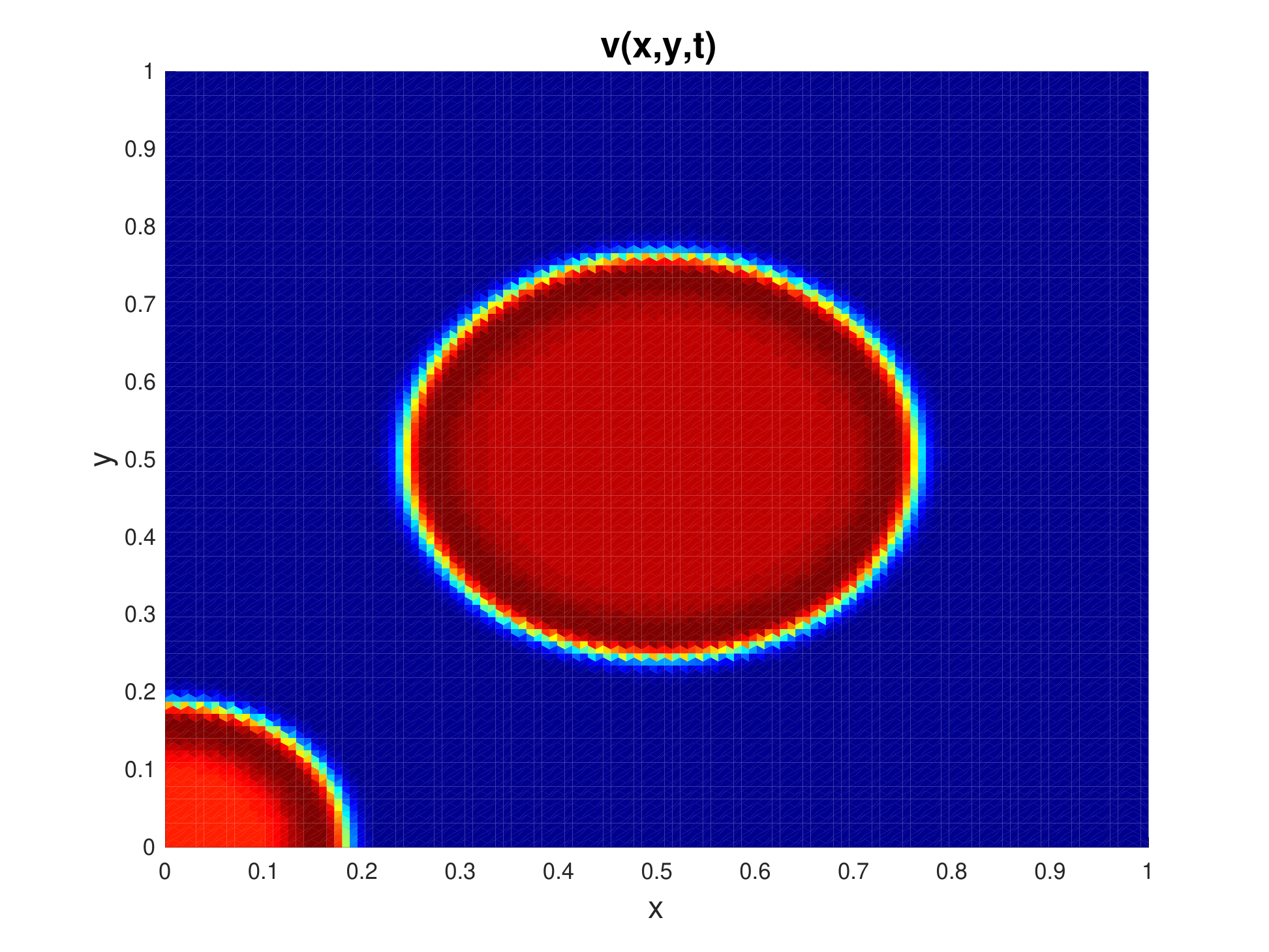} \\
$t=0.1$&$t=0.5$ &$t=1$
\\
\includegraphics[width=5cm, height=5cm]{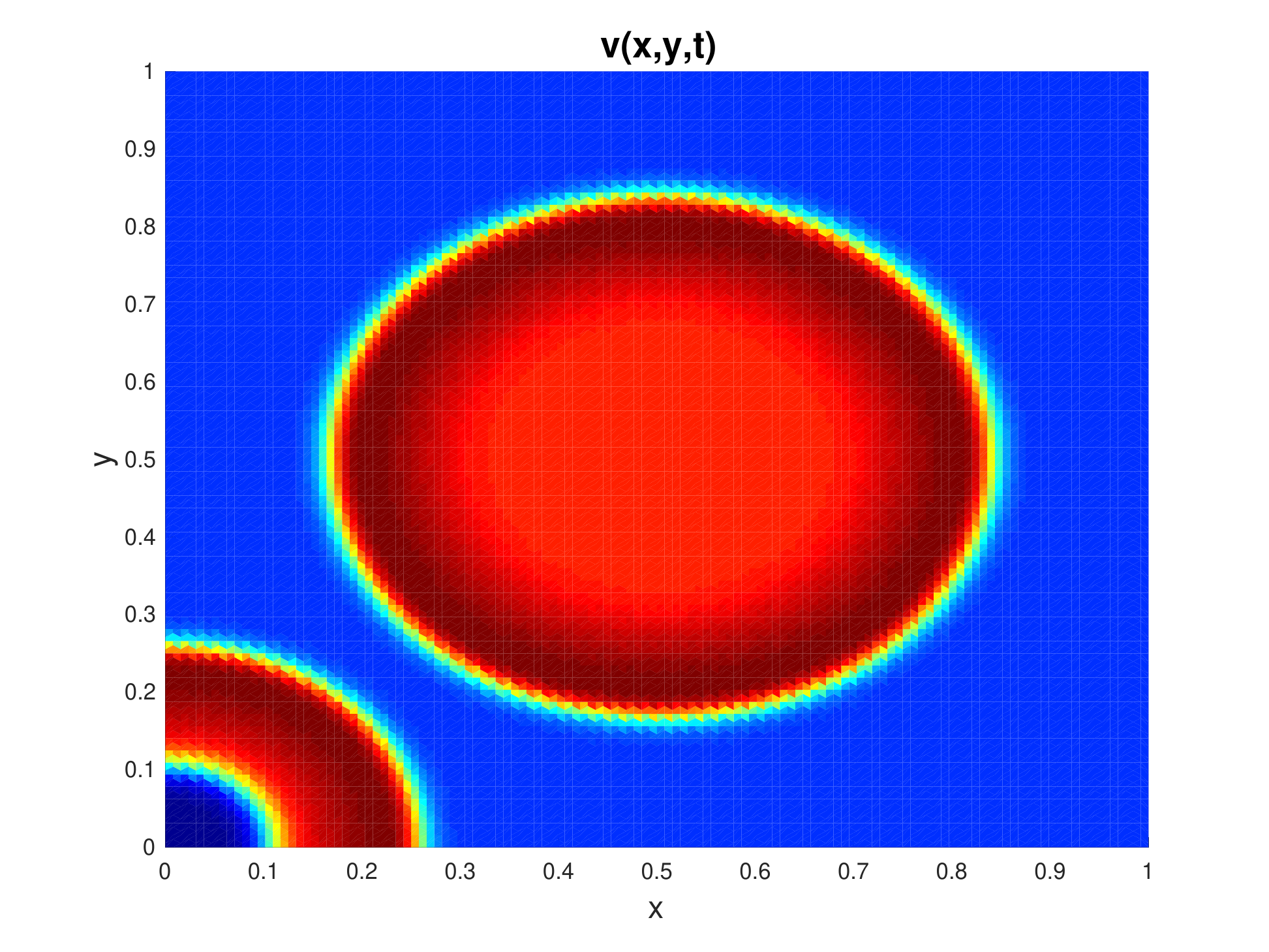} &
\includegraphics[width=5cm, height=5cm]{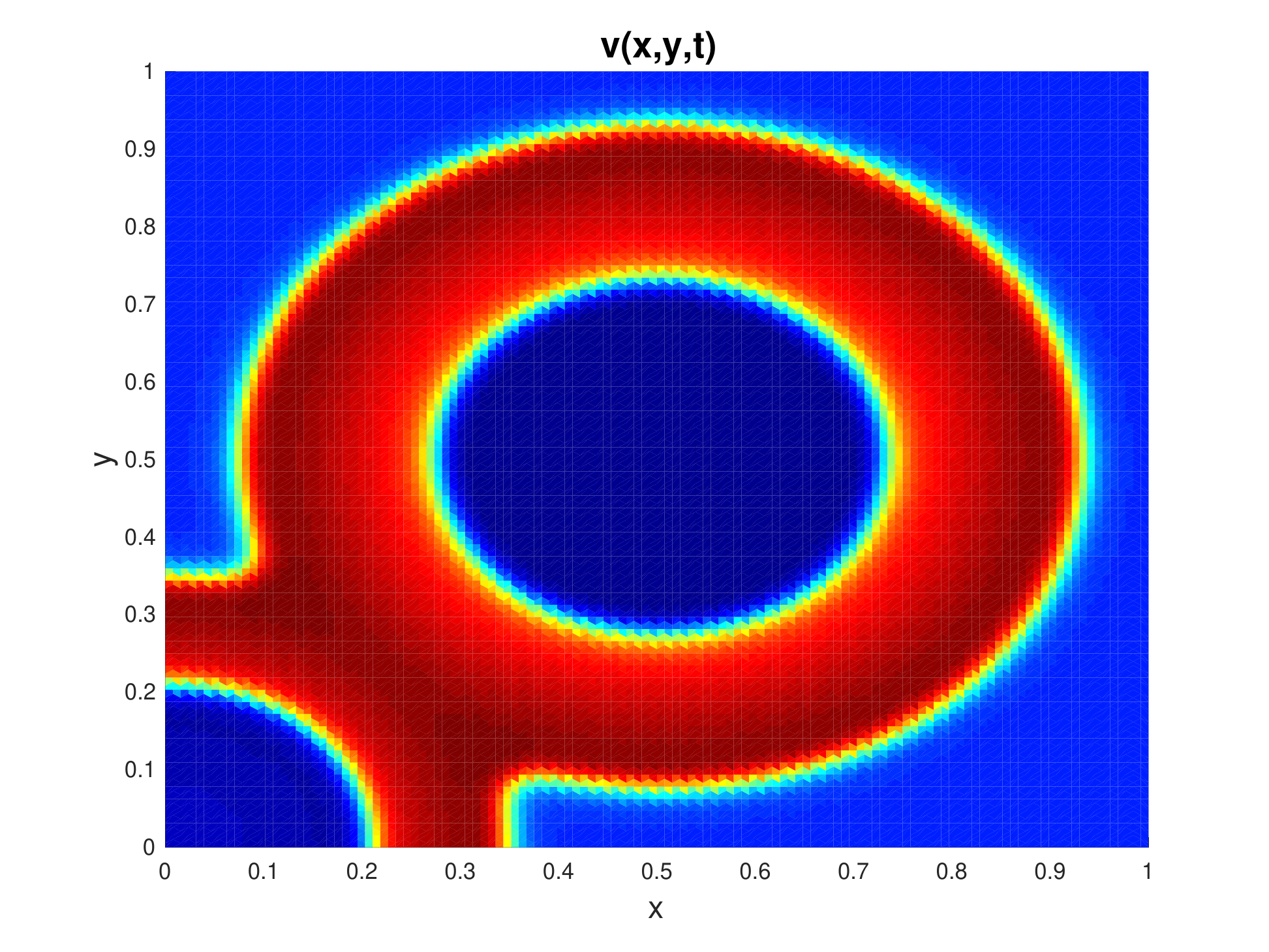} &
\includegraphics[width=5cm, height=5cm]{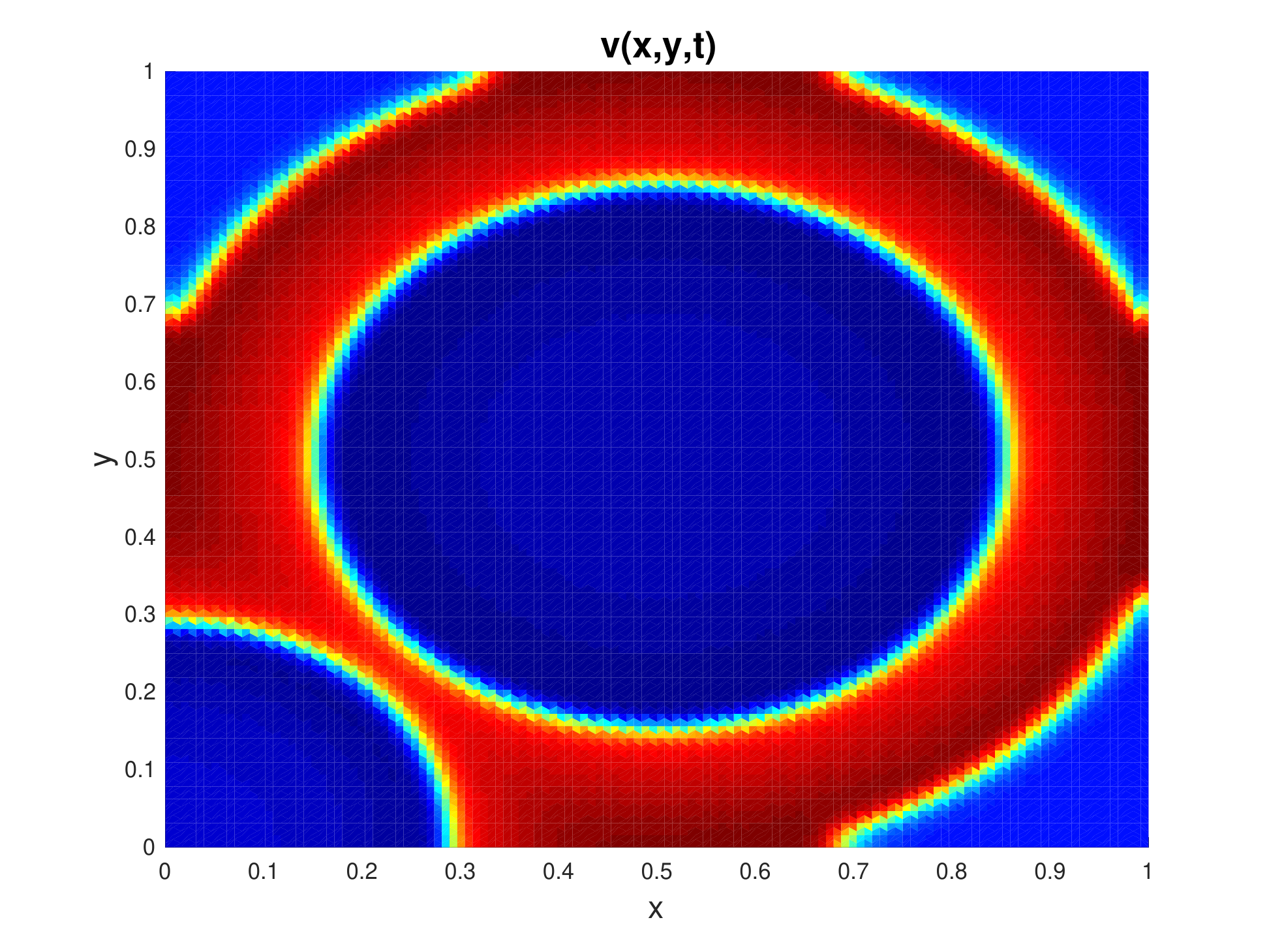} \\
$t=1.5$&$t=2$ &$t=2.5$
\\
\includegraphics[width=5cm, height=5cm]{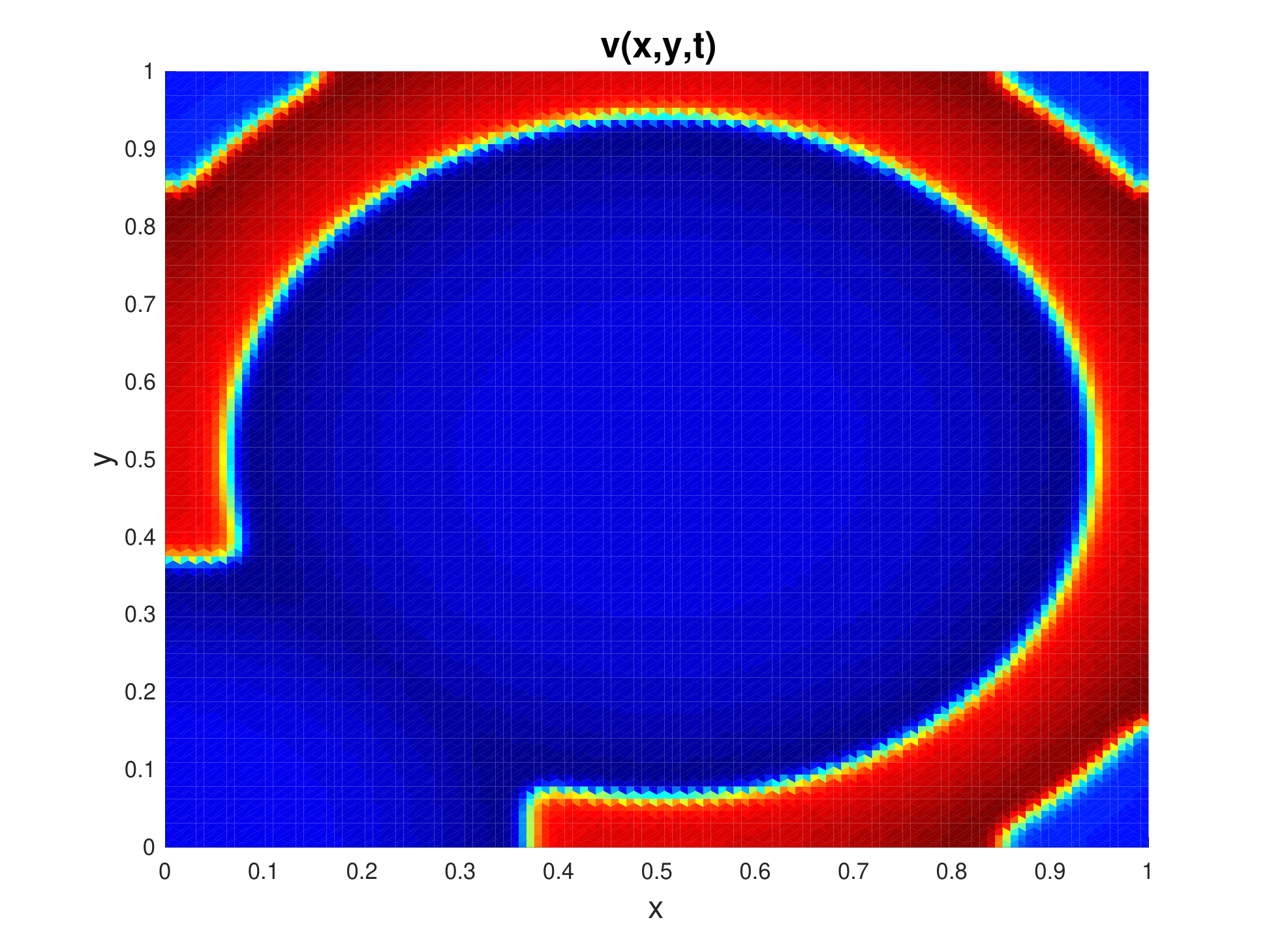} &
\includegraphics[width=5cm, height=5cm]{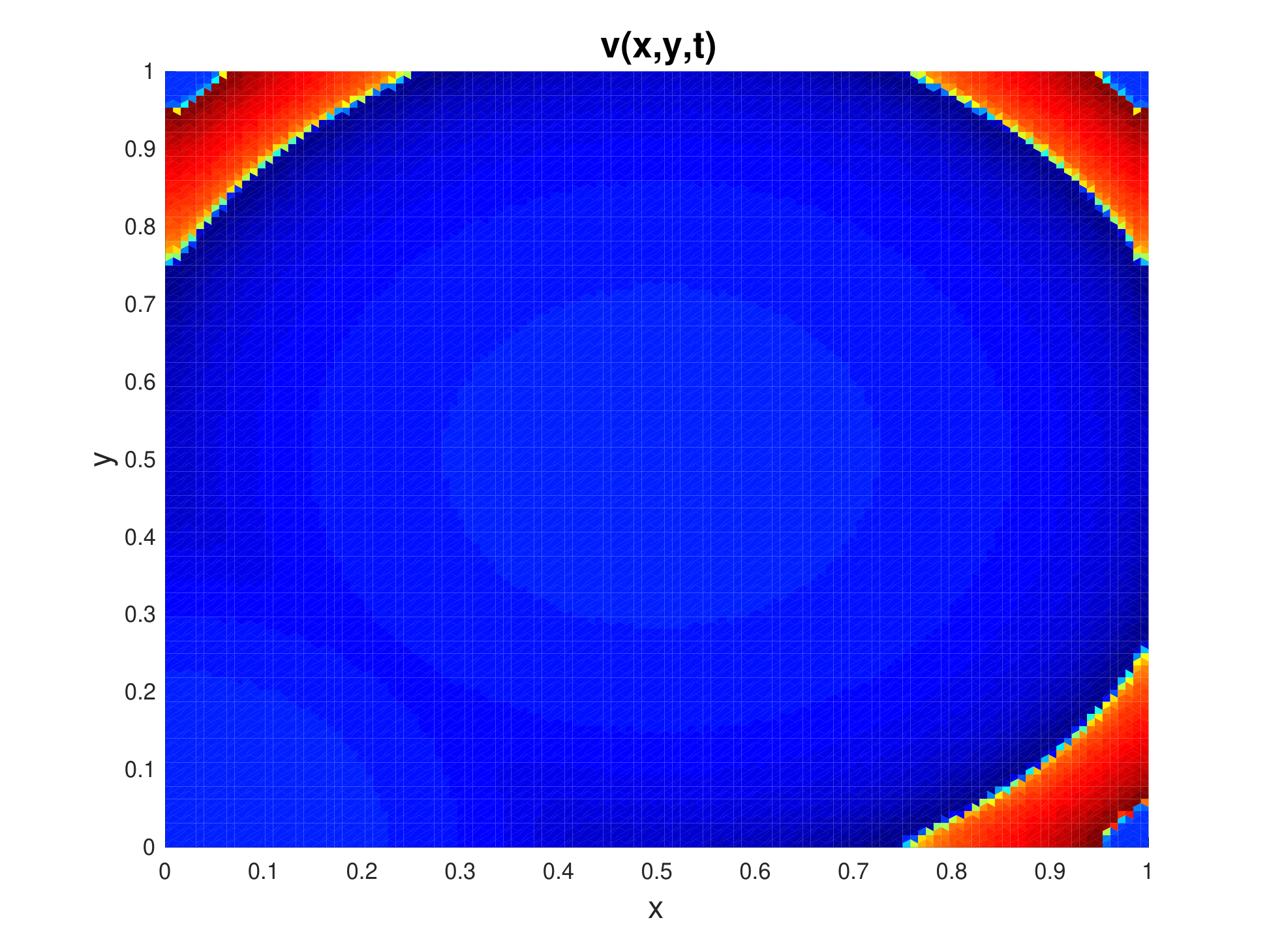} &
\includegraphics[width=5cm, height=5cm]{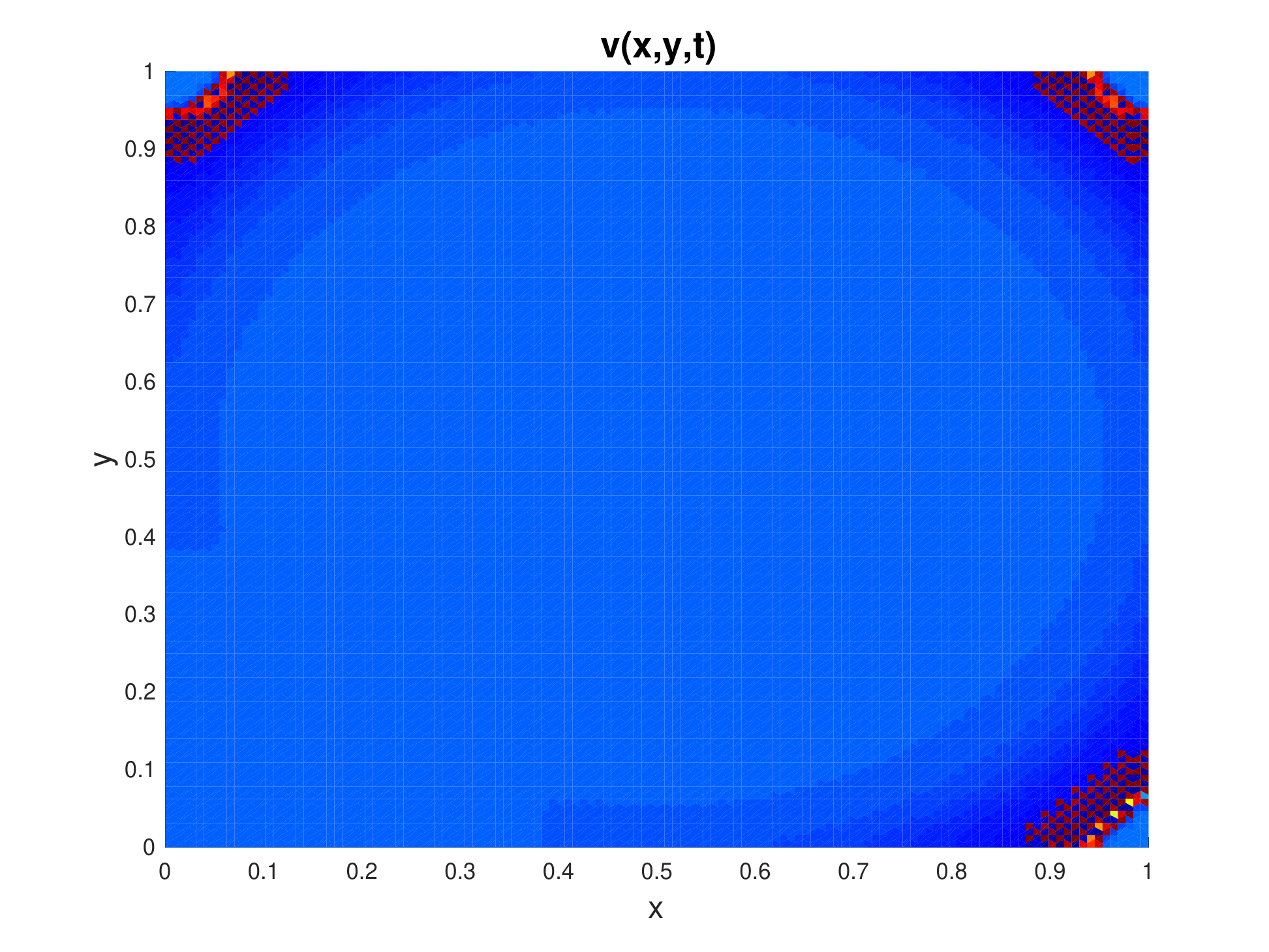}\\
$t=3$&$t=4$ &$t=5$ 
\end{tabular}
\caption{Numerical solution of the transmembrane potential $v$ for different times.}
\label{fig1}
\end{figure}

We show in Figure~\ref{fig1} the 
evolution of the transmembrane potential $v$ for different times.

\subsection{Example 3}

\begin{center}
\begin{figure}[ht]
\begin{tabular}{ccc}
\includegraphics[width=5cm, height=4.5cm]{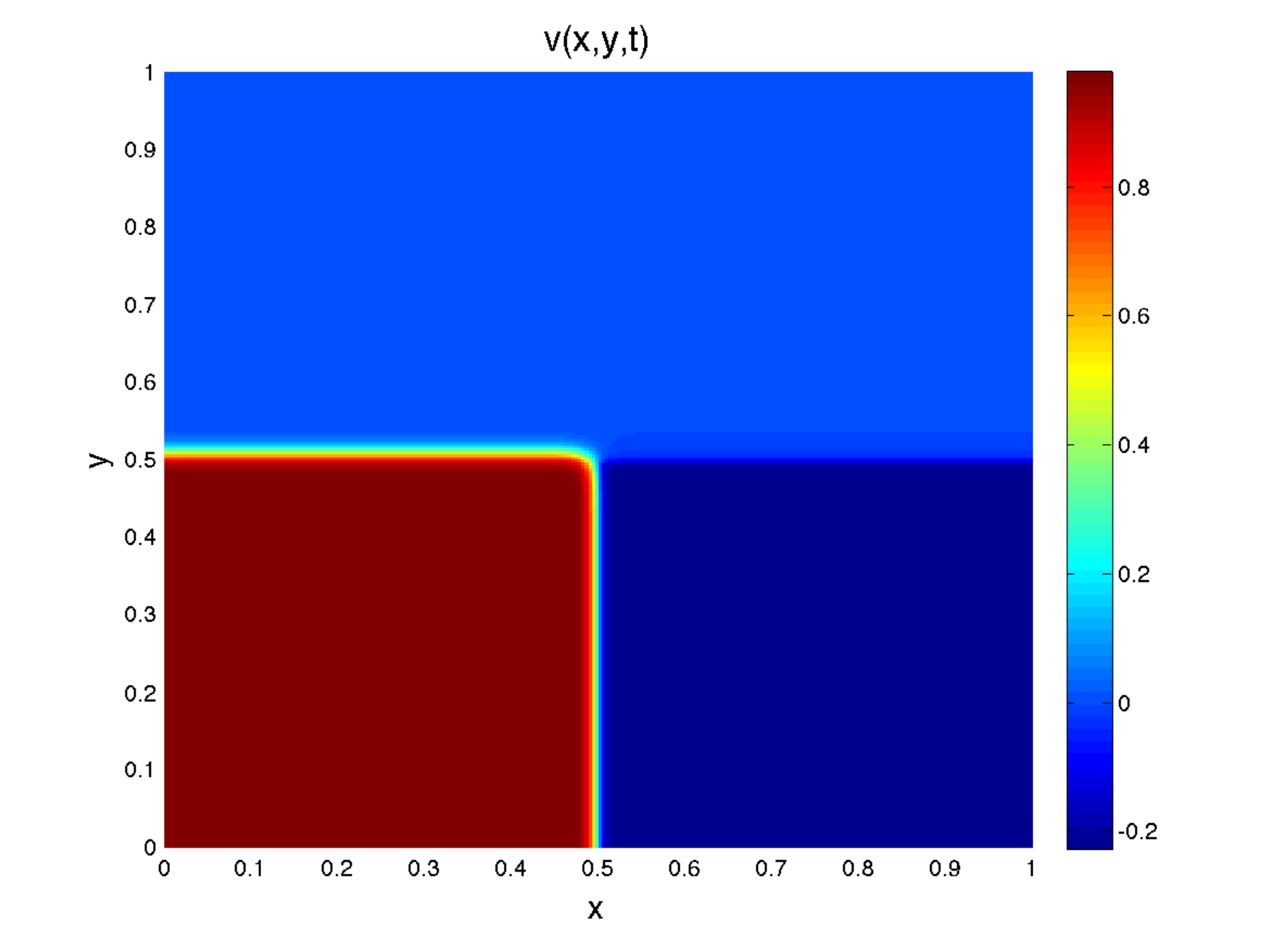} &
\includegraphics[width=5cm, height=4.5cm]{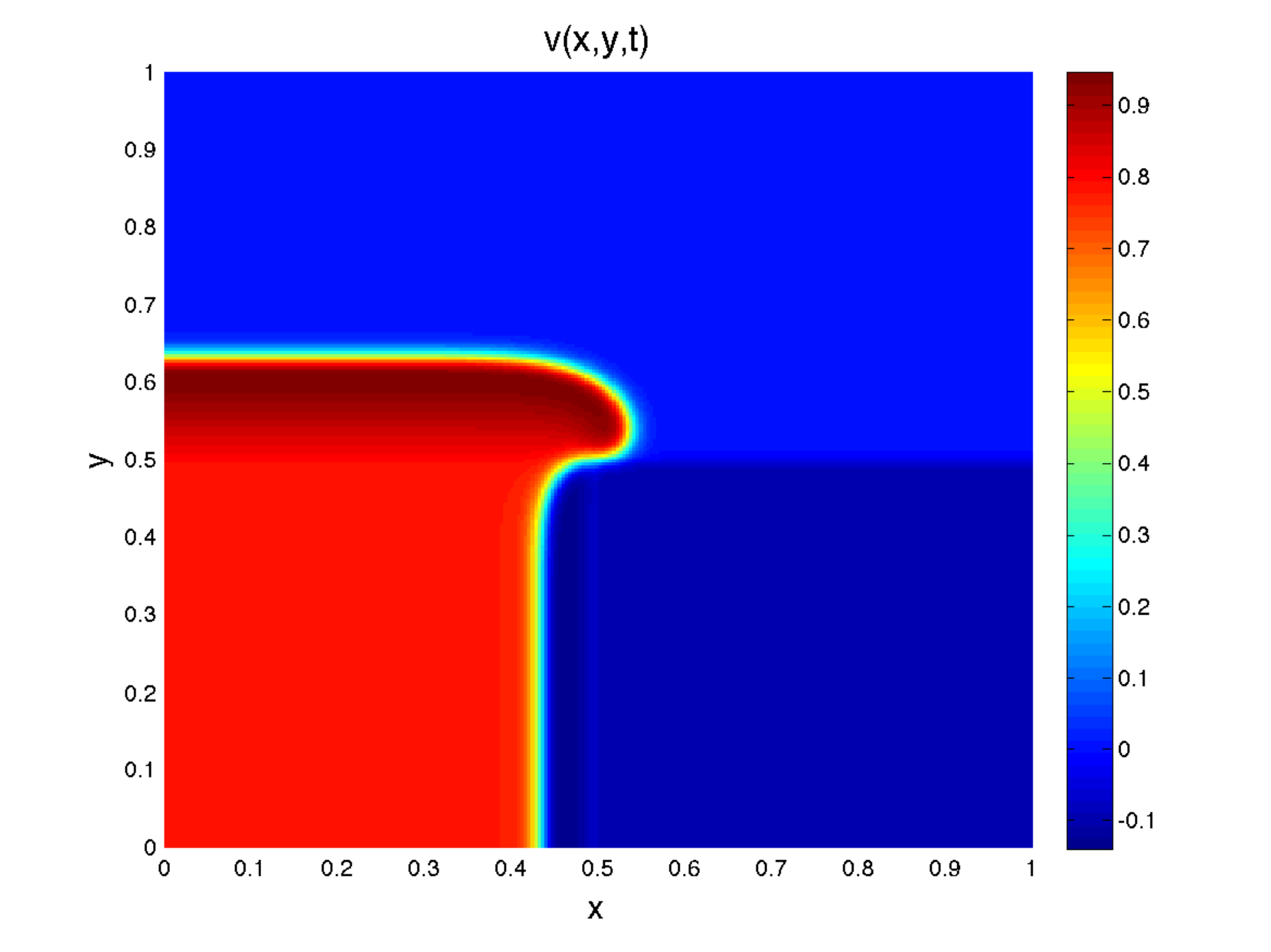} &
\includegraphics[width=5cm, height=4.5cm]{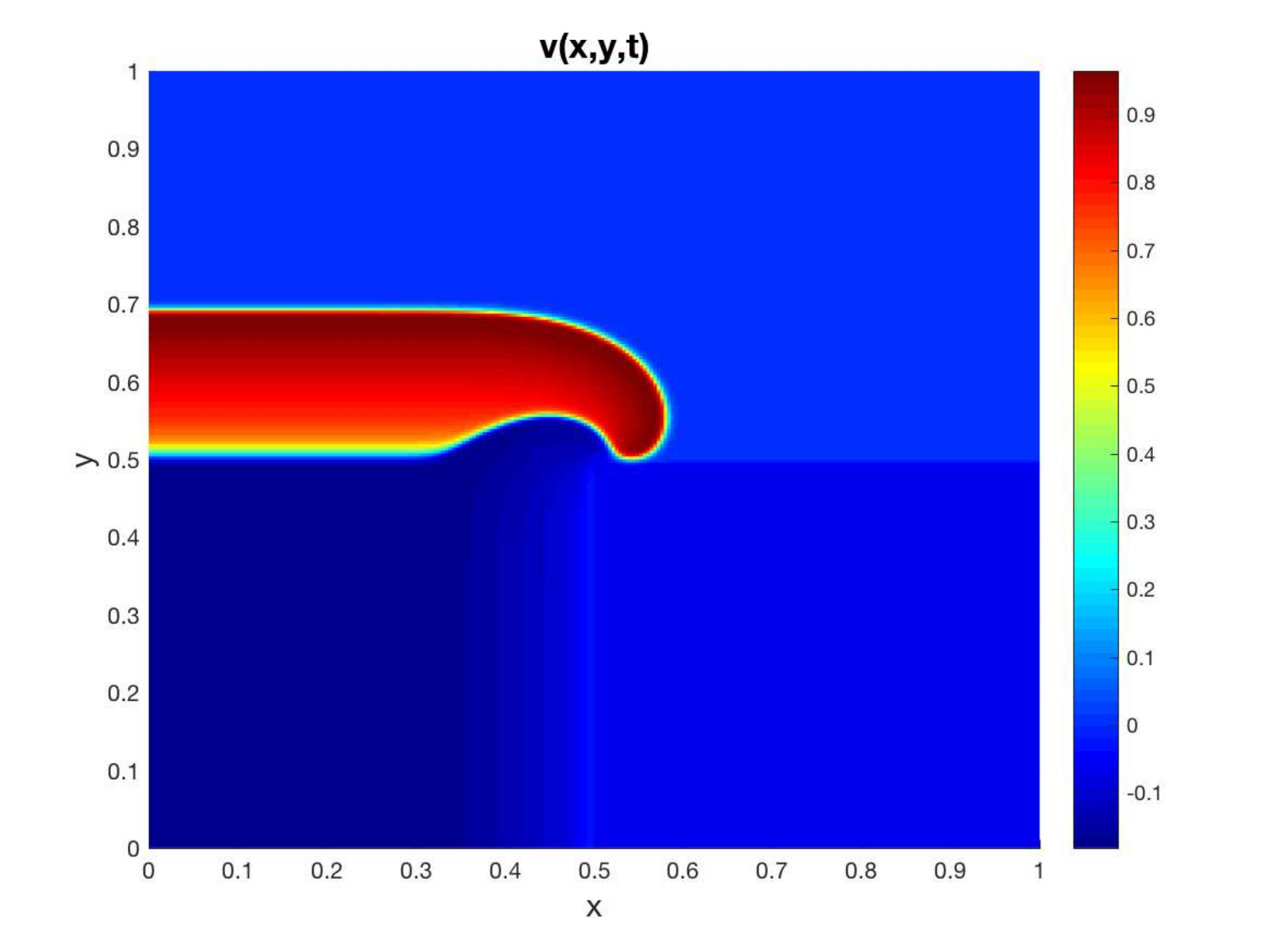}
\\
$t=0.1$ &$t=1.0$&$t=1.5$
\\
\includegraphics[width=5cm, height=4.5cm]{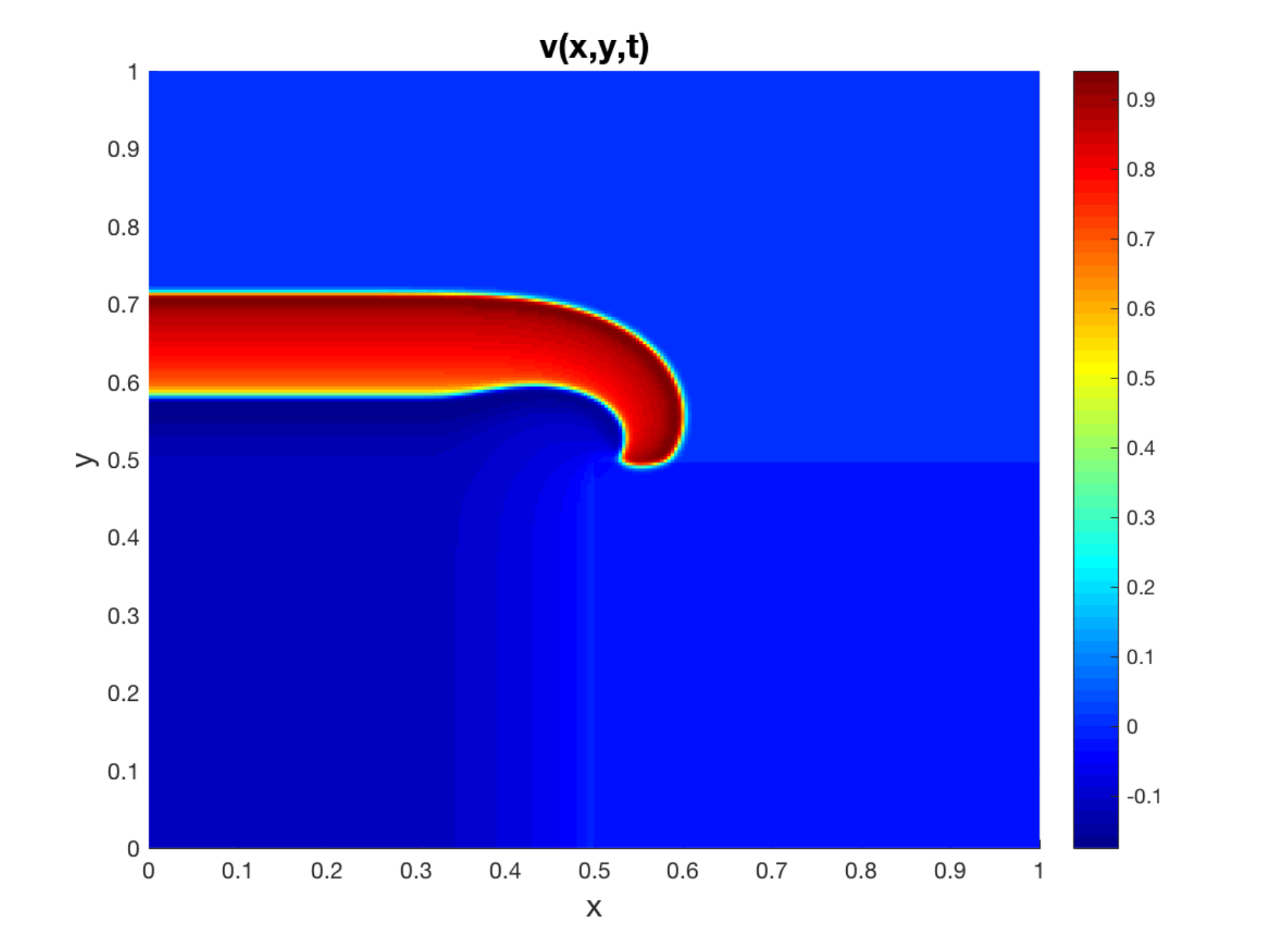} &
\includegraphics[width=5cm, height=4.5cm]{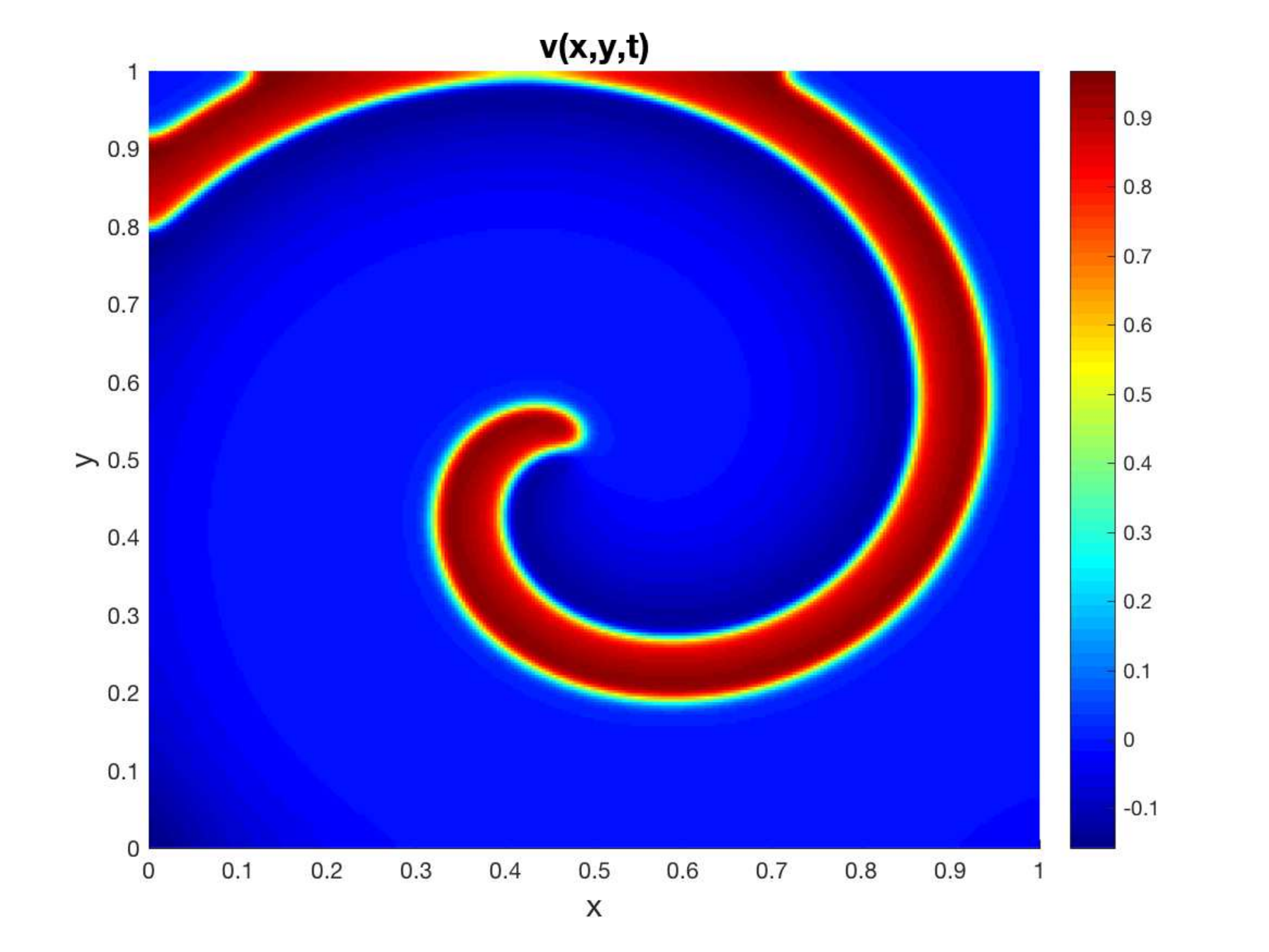} &
\includegraphics[width=5cm, height=4.5cm]{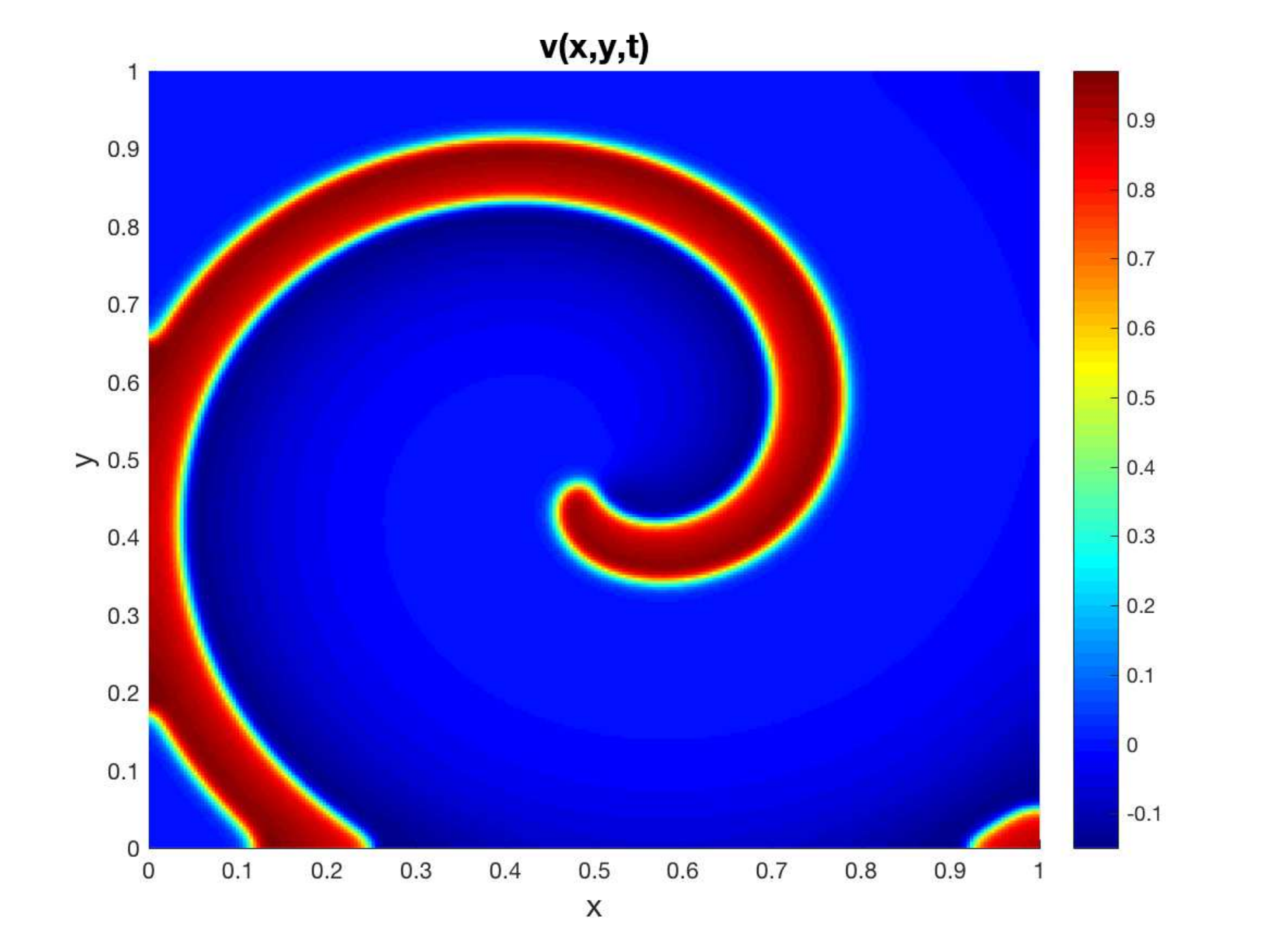} \\
\\
$t=2.0$ &$t=10.0$ &$t=15.0$
\end{tabular}
\caption{Numerical solution of the transmembrane potential $v$ for different times.}
\label{fig2}
\end{figure}
\end{center}

The aim of this test is to obtain the well-known periodic
spiral wave (see Figure~\ref{fig2}).
For this example, we use meshes $\mathcal T_h^3$ (with $h=1/128$)
on the domain $\Omega:=(0,1)^2$, and time interval $[0,15]$
(with $\Delta t=1/200$).
We will take the model constants as follow
$a=0.16875, b=1, \lambda=-100, \theta=0.25$.
Moreover, we consider the following initial data:

\begin{equation*}
v^0(x,y)=\left\{ \begin{array}{lcccc}
             1.4 &   if  & x < 0.5 & and & y< 0.5  \\
             \\ 0 &  otherwise, \\
             \end{array}
   \right. \quad
\end{equation*}
\begin{equation*}
   w^0(x,y)=\left\{ \begin{array}{lcccc}
             0.15 &   if  & x > 0.5 & and & y< 0.5  \\
             \\ 0 &  otherwise. \\
             \end{array}
   \right.
\end{equation*}

As expected the initial data
evolves to a spiral wave;
see Figure~\ref{fig2}. 

\section*{Acknowledgment}
V. Anaya was partially supported by
CONICYT-Chile through FONDECYT project 11160706;
D. Mora was partially supported by
CONICYT-Chile through FONDECYT project 1180913;
M. Sep\'ulveda was partially supported by
CONICYT-Chile through FONDECYT project 1180868,
and Basal, CMM, Universidad de Chile and
CI\textsuperscript{2}MA, Universidad de Concepci\'on

\bibliographystyle{amsplain}

\end{document}